\documentclass[number,citesort,dvips]{arxbj}
\usepackage{graphicx}

\aid{0}
\volume{17}
\issue{1}
\pubyear{2011}
\firstpage{60}
\lastpage{87}
\doi{10.3150/10-BEJ260}

\makeatletter
\DeclareMathAlphabet\mathcaligr{OMS}{cmsy}{m}{n}
\renewcommand{\mathcal}{\mathcaligr}
\newtheorem{lemma}{Lemma}[section]
\newtheorem{theorem}{Theorem}[section]
\makeatother

\begin{document}
\begin{frontmatter}

\title{Nonparametric regression with filtered data}
\runtitle{Nonparametric regression with filtered data}

\begin{aug}
\author[a]{\fnms{Oliver} \snm{Linton}\thanksref{a}\ead[label=e1]{o.linton@lse.ac.uk}},
\author[b]{\fnms{Enno} \snm{Mammen}\thanksref{b}\ead[label=e2]{emammen@rumms.uni-mannheim.de}},
\author[c]{\fnms{Jens Perch} \snm{Nielsen}\thanksref{c}\ead[label=e3]{Jens.Nielsen.1@city.ac.uk}}
\and\\
\author[d]{\fnms{Ingrid}~\snm{Van~Keilegom}\corref{}\thanksref{d}\ead[label=e4]{ingrid.vankeilegom@uclouvain.be}}
\runauthor{Linton, Mammen, Nielsen and Van Keilegom}
\address[a]{Department of Economics, London School of Economics,
Houghton Street, London
WC2A 2AE, UK.
\printead{e1}}
\address[b]{Department of Economics, University of Mannheim, L7, 3-5,
68131 Mannheim,
Germany.\\
\printead{e2}}
\address[c]{Cass Business School, 106, Bunhill Row, London EC1Y 8TZ,
UK.\\
\printead{e3}}
\address[d]{Institute of Statistics, Universit\'{e} catholique de
Louvain, Voie du Roman
Pays 20, B 1348 Louvain-la-Neuve, Belgium.
\printead{e4}}
\end{aug}

% HISTORY:
\received{\smonth{9} \syear{2008}}
\revised{\smonth{1} \syear{2010}}

% ABSTRACT
%
\begin{abstract}
We present a general principle for estimating a regression function
nonparametrically, allowing for a wide variety of data filtering, for example,
repeated left truncation and right censoring. Both the mean and the median
regression cases are considered. The method works by first estimating the
conditional hazard function or conditional survivor function and then
integrating. We also investigate improved methods that take account of model
structure such as independent errors and show that such methods can improve
performance when the model structure is true. We establish the pointwise
asymptotic normality of our estimators.
\end{abstract}

% KEYWORDS
%
\begin{keyword}
\kwd{censoring}
\kwd{counting process theory}
\kwd{hazard functions}
\kwd{kernel estimation}
\kwd{local linear estimation}
\kwd{truncation}
\end{keyword}

\end{frontmatter}

%s1 ###
\section{Introduction}

This paper concerns the nonparametric estimation of a regression
function $g(x)$
that regresses $Y$ on $X=x$, where the nonnegative variable $Y$ is subject
to various filtering schemes and where $X$ is an observed vector of
regressors. We consider both the mean and the median regression case. A
common particular case is the standard censored regression model $%
Y=g(X)+\varepsilon$, where $X$ is an observed $d$-dimensional vector of
regressors, $Y$ is subject to random right censoring and $\varepsilon$ is
an unobserved error satisfying $E(\varepsilon|X)=0$. We make two
contributions. First, we present a completely nonparametric estimation
methodology. This is done under more general censoring patterns than in
previous papers. Second, we assume that the error is independent of the
covariate and we show how to construct a more efficient estimator that
takes account of the common shape.

Parametric and semiparametric estimators of censored regression models
include Heckman \cite{Heckman1976}, Buckley and James  \cite{Buckley1979},
Koul, Susarla and Van Ryzin  \cite{Koul1981}, Powell \cite{Powell1984,Powell1986a,Powell1986b}, Duncan  \cite{Duncan1986}, Fernandez  \cite{Fernandez1986},
Horowitz \cite{Horowitz1986,Horowitz1988}, Ritov  \cite{Ritov1990}, Honor\'{e} and Powell \cite{Honor1994},
%Lewbel (1998a,b),
Buchinsky and Hahn  \cite{Buchinsky1998} and
Heuchenne and Van Keilegom  \cite{Heuchenne2007a}. Many of
these authors either assume $g(x)=\beta^{\top}x$ or some other parametric
form, provide estimates of average derivatives only up to an unknown
scale or assume that the error distribution is parametric. The fully
nonparametric $g(x)$ model we consider is important because of the
sensitivity of the parametric and semiparametric estimators to
misspecification of functional form. A small number of estimators exist for
nonparametric censored regression models, in most cases focusing on the
standard random censoring model. %Fan and Gijbels (1994) proposed a
%nonparametric censored regression estimator based on a local version of
%Buckley and James (1979). This estimator is consistent when the
%censoring
%point is drawn from a continuous distribution. Other possible
Dabrowska  \cite{Dabrowska1992} and Van Keilegom and Veraverbeke \cite{Keilegom1998} proposed
nonparametric censored regression estimators based on quantile methods.
Lewbel and Linton  \cite{Lewbel2002} considered the above standard censoring model,
except that the censoring time $C$ is taken to be a degenerate random
variable (i.e., it
is constant), while Heuchenne and Van Keilegom \cite{Heuchenne2007b,Heuchenne2008} considered the
standard model when it is supposed that $\varepsilon$ is independent
of $X$.

In this paper, we propose a unified approach to the estimation of the
regression function from filtered data. Filtering, for example, left
truncation or right censoring, means that even though some information is
available about $Y$, $Y$ itself is sometimes not observed, even though
$X$ is
observed. It is imperative for us that our estimation principles are natural
and well known in the simple case of independent identically distributed
errors with no filtering. Our approach makes use of tools from the
field of
counting process theory; see
   \cite{Andersen1992} and
   \cite{Fleming1991}.

First, we recognize that the generic regression model can be reformulated
through the counting process $N(y)=I(Y<y)$ such that $Y=\int
_{0}^{\infty
}I(Y>y)\,\mathrm{d}y=\int_{0}^{\infty}yN(\mathrm{d}y).$ The advantage of the counting process
approach is that it readily lends itself to quite general filtering
mechanisms, allowing for complicated left truncation and right censoring
patterns.

We reformulate the regression model in terms of a counting process $N$
having stochastic intensity function
\[
\lambda(y)=\alpha_{X}(y)Z(y)
\]
with respect to the increasing, right-continuous and complete
filtration $%
\mathcal{F}_{y}= \{ X,N(u)\mid 0<u\leq y \}$. Here, $Z(y)=1-N(y)$ and $%
\alpha_x(y)$ is the conditional hazard function of $Y$ given that $X=x$.
%} ) , \label{eq:multi}
%where
%_{0}|x}(u)/(1-\int_{0}^{u}f_{\varepsilon_{0}|x}(s)ds)
%is the conditional hazard function of the error term $\varepsilon
%_{0}=Y/E(Y|X)$.
With these definitions, we have that the conditional mean is given by
%
%e1 ###
\begin{equation}
g_{\mathrm{mn}}(x) = E(Y|X=x) = - \int_0^\infty y S_x(\mathrm{d}y) = \int_{0}^{\infty
}u\alpha_{x}(u)\exp\biggl(-\int_{0}^{u}\alpha_{x}(v)\,\mathrm{d}v \biggr) \,\mathrm{d}u
\label{eq:mean}
\end{equation}
and the conditional median is given by
%
%e2 ###
\begin{equation}
g_{\mathrm{med}}(x) = S_x^{-1}(0.5), \label{eq:median}
\end{equation}
where the relation between the conditional survival function $S_x(\cdot)$
and the conditional hazard function $\alpha_x(\cdot)$ is given by
\[
S_x(y) = \exp\biggl\{ - \int_0^y \alpha_x(u) \,\mathrm{d}u \biggr\}.
\]
%
%Recall that the conditional survival function satisfies two well known
%relations
%S_{x}(y) &=&\underset{w\leq y}{\prod}\{1-\Lambda_{x}(dw)\}
%&& \notag\\
%&=&\exp\{-\Lambda_{x}(y)\}, \label{s2}
%where $\prod$ denotes the product integral [Gill (1980)] and $\Lambda
%_{x}(y)=\int_{0}^{y}\alpha_{x}(u)du$ is the integrated conditional
%hazard
%function$.$ Therefore, we can write
%g(x)=-\int_{0}^{\infty}yS_{x}(\mathrm{d}y), \label{ds}
%where $S_{x}(y)$ is defined in terms of $\alpha_{x}(y)$ by (\ref{s1})
%or by
%(\ref{s2}).
This connection between the hazard function and the regression function is
the basis of our estimation.

For the first contribution of this paper, we consider $\alpha_{x}(y)$ as
estimated from a local constant least-squares principle or a local linear
least-squares principle. Plugging these estimators into the expressions
(\ref%
{eq:mean}) and (\ref{eq:median}) results in, respectively, a local
constant $\hat{g}_{C}$ and
a local linear estimator $\hat{g}_{L}$ of the
conditional mean or median.
%We note that $\widehat{g}_{C}=\widehat{m}_{C}$ and $\widehat{g}_{L}=%
%is
%used where we do not smooth in the $y$-direction. This is a quite
%powerful
%statement since it implies that
It is important to note that in the absence of filtering, the traditional
local constant and local linear kernel regression estimators are special
cases of the estimators $\hat{g}_{C}$ and $\hat{g}_{L}$.
%local constant
%and local linear kernel hazard estimators of this paper.

The second contribution of this paper is concerned with the estimation of
the functions $g_{\mathrm{mn}}(\cdot)$ and $g_{\mathrm{med}}(\cdot)$ when some
structure is
imposed on the model.
%We also consider the case where the error distribution is generated by
%the
%same underlying shape (and independent of the covariates).
If there is a substantial level of filtering, then one can envision areas
where truncation or censoring imply that we do not have local
information on
the entire shape of the error distribution around every~$x$. One can
alleviate this by imposing assumptions on the shape of these local error
distributions. The simplest model assumption in this connection is the
multiplicative regression model
%simply that $\alpha_{\varepsilon_{0}|x}$ does not
%depend on $x$. This is
%
%e3 ###
\begin{equation}
Y=g(X)\varepsilon_{0}, \label{eq:genmodel2}
\end{equation}
where the error term $\varepsilon_{0}$ is independent of $X$ and has mean
or median equal to one, and where $g(X)$ is either $g_{\mathrm{mn}}(X)$ or $%
g_{\mathrm{med}}(X) $. Under this model,
%
%e4 ###
\begin{equation}
\alpha_{\varepsilon_{0}|x}\equiv\alpha_{0} \label{eq:alfanul}
\end{equation}
for some function $\alpha_{0},$ where $\alpha_{\varepsilon_0|x}$ is the
conditional hazard function of $\varepsilon_0$ given that $X=x$.
%This is equivalent to a multiplicative regression model

If model (\ref{eq:genmodel2}) is true, then it can be used to improve
estimation, even in the case without filtering; see
   \cite{Tibshirani1984}. Our
estimation strategy in this case is sequential. We first obtain the
unrestricted estimator $\hat{g}(\cdot)=\hat{g}_{\mathrm{mn}}(\cdot)$ or $%
\hat{g}_{\mathrm{med}}(\cdot)$ described above. We then use the relation
%
%e5 ###
\begin{equation}
\alpha_{x}(y)=\frac{1}{g(x)}\alpha_{0} \biggl( \frac{y}{g(x)} \biggr)
\label{eq:multi}
\end{equation}
or, equivalently, $\alpha_{0}(u)=g(x)\alpha_{x}(ug(x))$ to obtain an
estimate for $\alpha_{0}(\cdot).$ We use a minimum chi-squared approach to do
this optimally, which involves replacing $g(x)$ by $\hat{g}(x)$ and
$%
\alpha_{x}(y)$ by the completely nonparametric estimator $\hat
{\alpha}%
_{x}(y)$. Given an estimator of $\alpha_{0}(\cdot)$, we then obtain a new
estimator of $g(x)$ using the minimum chi-squared approach, again based on
the relation $\alpha_{x}(y)=\alpha_{0} ( y/g(x) ) /g(x)$, but now
replacing $\alpha_{0}(u)$ by $\hat{\alpha}_{0}(u)$ and $\alpha_{x}(y)$
by $\hat{\alpha}_{x}(y)$. We will argue that our estimator
fulfills a
local efficiency criterion. Van Keilegom and Akritas \cite{Keilegom1999} and Heuchenne
and Van Keilegom \cite{Heuchenne2007b,Heuchenne2008} discuss estimation of $S_{x}(y)$ and $E(Y|X=x)
$, respectively, in the additive error model when $Y-E(Y|X)$ is
independent of
$X$. In the first two papers, $S_{x}(y)$ or $E(Y|X=x)$, respectively, is
written as a functional of the error distribution and of the
distribution of
the covariates. The estimator is based on plugging in estimates of these
distributions. In the last paper, censored observations are replaced by
synthetic data points. In all three of these papers, efficiency issues
are not
discussed and the analysis is restricted to the case of random right
censoring.

The outline of the paper is as follows. In Section~\ref{sec2}, we describe the
theoretical background in terms of the counting process formulation,
including the important special case of filtered data. In Section~\ref{sec3}, we
introduce our approach to regression based on filtered data in the general
situation, where we do not restrict the functional form of the error
distribution. We present the local constant case in detail; the local linear
case is given in the \hyperref[appm]{Appendix}. The more efficient estimator (at least when
the assumption is correct) based on the assumption on the functional form
(assumption (\ref{eq:alfanul})) is introduced in Section~\ref{sec4}, where we also
give its asymptotic distribution. In Section~\ref{sec5}, we present a small
simulation study. In the \hyperref[appm]{Appendix}, we give the proofs of the main
distribution results contained in the text.

%s2 ###
\section{The counting process framework}\label{sec2}

Let $(X_i,Y_i)$, $i=1,\ldots,n$, be $n$ i.i.d. replications of the random
vector $(X,Y)$, where the response $Y_i$ is subject to filtering and
therefore possibly unobserved, and the covariate $X_i=(X_{i1},\ldots,X_{id})$
is completely observed.

%s2.1 ###
\subsection{The unfiltered case}

Define $N_{i}(y)=I(Y_{i}<y)$ for all $y$ in the support of $Y_{i}$.
Then $%
\mathbf{N}=(N_{1},\ldots,N_{n})$ is an $n$-dimensional counting process
with respect to possibly different, increasing, right-continuous, complete
filtrations $\mathcal{F}_{y}^{i}$; see  \cite{Andersen1992},
page 60. We
assume that with respect to the filtration, $N_{i}$ has stochastic intensity
%
%e6 ###
\begin{equation}
\lambda_{i}(y)=\alpha_{X_{i}}(y)Z_{i}(y), \label{eq:inten}
\end{equation}
where $Z_{i}(y)=I(Y_{i}\geq y)$ is a predictable process taking values
in $%
\{0,1\}$.
%indicating (by the value $1$) when the $i$'th individual is under
%risk.
We have not restricted the conditional distribution of $S_{X_i}$ and the
functional form of the conditional hazard function $\alpha_{X_i}$ is
likewise unrestricted. With these definitions, $\lambda_i $ is predictable,
and the processes $M_{i}(y)=N_{i}(y)-\Lambda_{i}(y),$ $i=1,\ldots,n,$ and
compensators $\Lambda_{i}(y)=\int_{0}^{y}\lambda_{i}(s)\,\mathrm{d}s,$ are
square-integrable local martingales on the support of $Y_i$.

We can allow this extremely general model description since the martingale
central limit theorem dating back to
Rebolledo  \cite{Rebolledo1980} can be applied in this
context; see  \cite{Andersen1992}, pages 82--85. Our
framework is
sufficiently general to include a number of interdependencies,
including a variety of time series analyses.
%Therefore, even in the case without filtering, our model is
%more general than the model considered in Mammen and Nielsen (2007).

%s2.2 ###
\subsection{The filtered case}

In this section, we follow
Andersen  \cite{Andersen1988},
page 50.
Let $C_{i}(y)$ be a predictable process taking values in $\{0,1\}$,
indicating (by the value $1$) when the $i$th individual is at risk. Note
that the predictability condition of $C_{i}(y)$ allows it to depend on $
X_{i}=(X_{i1},\ldots,X_{id})$ in every possible way. Let
\[
\overline{N}_{i}(y)=\int_{0}^{y}C_{i}(s)\,\mathrm{d}N_{i}(s)
\]
be the filtered counting process and introduce the filtered filtration $
\overline{\mathcal{F}}_{y}=\sigma(\overline{N}(s),X,CZ(s);s\leq y).$
The random intensity process $\overline{\lambda}_i$ is then
\[
\overline{\lambda}_{i}(y)=\alpha_{X_i} (y)C_{i}(y)Z_{i}(y)
\]
and the integrated random intensity process is
\[
\overline{\Lambda}_{i}(y)=\int_{0}^{y}\overline{\lambda}%
_{i}(s)\,\mathrm{d}s=\int_{0}^{y}\alpha_{X_i}
(s)C_{i}(s)Z_{i}(s)\,\mathrm{d}s=\int_{0}^{y}C_{i}(s)\,\mathrm{d}\Lambda_{i}(s).
\]
With these definitions, $\overline{M}_{i}(y)=\overline
{N}_{i}(y)-\overline{%
\Lambda}_{i}(y)$ is a square-integrable martingale with respect to the
filtration $ ( \overline{\mathcal{F}}_{y} ) _{y\geq0}.$ Note that,
in the filtered case, $Z_i(y) = I(Y_i \ge y)$ is not always observed, but
the product $(C_iZ_i)(y)$ is always observable.

%s3 ###
\section{Estimation under the completely nonparametric model}\label{sec3}

In this section, local constant and local linear estimators under the general
nonparametric model are given. These estimators take the local constant and
the local linear marker-dependent kernel hazard estimators of
Nielsen and Linton  \cite{Nielsen1995} and
Nielsen  \cite{Nielsen1998} as their starting point. In the
special case of
no filtering, this results in the convenient property that the regression
estimator based on the local constant hazard estimator is the well-known
local constant regression estimator, the Nadaraya--Watson estimator,
and the
local linear hazard estimator results in the local linear regression
estimator; see, for example,
  \cite{Fan1996}.

Let $K$ be a $d$-dimensional kernel, $k$ be a one-dimensional kernel, $%
b=(b_{1},\ldots,b_{d})$ be a $d$-dimensional bandwidth vector and $h$
be a
one-dimensional bandwidth. For any real $u$ and any $d$-dimensional
vector $%
x=(x_1,\ldots,x_d)$, define $k_{h}(u)=k(u/h)/h$ and $%
K_{b}(x)=|b|^{-1}K(x/b), $ where $x/b=(x_{1}/b_{1},\ldots
,x_{d}/b_{d})$ and
$|b|=\prod_{j=1}^{d}b_{j} $. The estimator suggested by Nielsen and Linton
(1995) is
%
%e7 ###
\begin{equation}
\hat{\alpha}_{x,C}(y)=\frac{O_{x,y}^{C}}{E_{x,y}^{C}},
\label{eq:localconstant}
\end{equation}
where
\begin{eqnarray*}
O_{x,y}^{C}&=&n^{-1}\sum_{i=1}^{n}\int K_{b} ( x-X_{i} ) k_{h}(y-u)\,\mathrm{d}%
\overline{N}_{i}(u),
\\
E_{x,y}^{C}&=&n^{-1}\sum_{i=1}^{n}\int K_{b} ( x-X_{i} )
k_{h}(y-u)C_{i}(u)Z_{i}(u)\,\mathrm{d}u.
\end{eqnarray*}
This estimator was identified as a local constant least-squares
estimator in
   \cite{Nielsen1998}. The super/subscript $C$ stands for local constant smoothing.
Below, we will also introduce estimators based on local linear smoothing.
This will be indicated by a super/subscript $L$ in the notation.

We wish to estimate the conditional integrated hazard $A
_{x}(y)=\int_{0}^{y}\alpha_x(u)\,\mathrm{d}u.$ We could just integrate $\hat
{\alpha
}_{x,C}(y)$ with respect to $y,$ but a better strategy is to first let the
bandwidth $h\rightarrow0,$ which eliminates redundant smoothing$.$ The
resulting estimator is
%e8 ###
\begin{equation} \label{def1S}
\widehat{A}_{x,C}(y)=\lim_{h\rightarrow0} \int_{0}^{y}\hat{%
\alpha}_{x,C}(u)\,\mathrm{d}u=\sum_{i=1}^{n}K_{b} ( x-X_{i} ) \int_{0}^{y}%
\frac{\mathrm{d}\overline{N}_{i}(u)}{\sum_{j=1}^{n}K_{b}(x-X_{j})C_{j}(u)Z_{j}(u)}.
\end{equation}
Note that $\widehat{A}_{x,C}(y)$ equals the estimator of $A _{x}(y)$
proposed by
Beran  \cite{Beran1981} and Dabrowska  \cite{Dabrowska1987} in the case of random
censoring. We then estimate the conditional survivor function
$S_{x}(y)$ by
the product limit estimator of %Johansen (1987) - see also
Johansen and Gill  \cite{Johansen1990}; see   \cite{Andersen1992},
that is,
%
%e9 ###
\begin{equation} \label{def2S}
\widehat{S}_{x,C}(y)=\prod_{0\leq w\leq y} \{ 1-\widehat{A }%
_{x,C}(\mathrm{d}w) \}
\end{equation}
for $y \le T$, where $T$ satisfies assumption (A) below. The local constant
estimator of $g_{\mathrm{mn}}^T(x)=E(Y I(Y \le T)|X=x)$ is%
%
%e10 ###
\begin{equation}
\hat{g}_{C,\mathrm{mn}}^T(x)=-\int_{0}^{T}y\widehat{S}_{x,C}(\mathrm{d}y). \label{gT}
\end{equation}
A local constant estimator of $g_{\mathrm{med}}(x)=\operatorname{med}(Y|X=x)$ is given by
\[
\hat{g}_{C,\mathrm{med}}(x)=\widehat{S}_{x,C}^{-1}(0.5),
\]
where for any $0<p<1$, $\widehat{S}_{x,C}^{-1}(p)=\inf\{y\dvtx \widehat{S}
_{x,C}(y)\leq1-p\}$.

Another option would have been to define $\overline{S}_{x}(y)=\exp\{-%
\widehat{A}_{x,C}(y)\}$ in the above formula. The advantage of the weighted
product limit estimator is that we arrive at exactly the extension of the
Kaplan--Meier estimator to filtered data in the absence of covariates
and at
the weighted empirical distribution function  \cite{Stone1977}  in the absence
of filtering. As a consequence, (\ref{gT}) reduces to the well-known
Nadaraya--Watson estimator when $T=\infty$ and when all data are completely
observed.

In a similar way, the local linear estimators of $S_x(y)$,
$g_{\mathrm{mn}}^T(x)$ and $%
g_{\mathrm{med}}(x)$, denoted $\widehat{S}_{x,L}(y)$, $\hat{g}_{L,\mathrm{mn}}^T(x)$
and $%
\hat{g}_{L,\mathrm{med}}(x)$, respectively, can be defined. We refer to the
\hyperref[appm]{Appendix} for their precise definitions.

For the asymptotic properties of the unrestricted estimators $\hat{g}
_{C,\mathrm{mn}}^T(x)$ and $\hat{g}_{L,\mathrm{mn}}^T(x)$ of $g_{\mathrm{mn}}^T(x)$, we need to
assume the following for $x \in R_X$, where $R_X$ is a bounded interval in
the interior of the support of $X$. All of our results are stated for the
special case of a one-dimensional covariate~$X$, $d=1$. The results can be
easily generalized to a multivariate setting.

\begin{enumerate}[(D6)]
\item[(D1)] The derivatives ${\frac{\partial^2\alpha_x(u)
}{\partial
x^2}} $
and ${\frac{\partial\alpha_x(u) }{\partial x}} $ exist and are
uniformly continuous in $x\in R_X,u\in\lbrack0,T]$.
\item[(D2)] The kernel $K$ is symmetric, continuous and has bounded support.
The bandwidth $b$ satisfies $b\rightarrow0$, $nb \rightarrow\infty$
and $%
nb^{5}=\mathrm{O}(1)$.
\item[(D3)] The truncation variable $T$ is such that $\inf_{x\in
R_X,u\in
\lbrack0,T]}\varphi_{x}(u)>0$.
\item[(D4)] There exists a continuous function $\varphi_{x}(y)$ such that
\begin{eqnarray*}
\sup_{y\in\lbrack0,T]} \Biggl\vert\frac{1}{n}%
\sum_{i=1}^{n}K_{b} ( x-X_{i} ) C_{i}(y)Z_{i}(y)-\varphi
_{x}(y)\Biggr \vert&\stackrel{P}{\rightarrow}&0, \\
\sup_{y\in\lbrack0,T]} \Biggl\vert\frac{1}{n}\sum_{i=1}^{n}%
\frac{ ( x-X_{i} )^2} {b^{2}} K_{b} ( x-X_{i} )
C_{i}(y)Z_{i}(y)-\frac{1} {2}\mu_{2}(K) \varphi_{x}(y) \Biggr\vert
&\stackrel{P}{\rightarrow}&0,
\end{eqnarray*}
where $\mu_{2}(K)= \int u^2 K(u) \,\mathrm{d}u$.
\item[(D5)] The derivative $\frac{\partial\varphi_{x}(y)}{\partial x}$
exists and is continuous. It holds that
\[
\sup_{y\in\lbrack0,T]} \Biggl\vert\frac{1}{n}\sum_{i=1}^{n}%
( X_{i}-x ) b^{-2} K_{b} ( x-X_{i} ) C_{i}(y)Z_{i}(y)-\mu
_{2}(K) \frac{\partial\varphi_{x}(y)}{\partial x} \Biggr\vert\stackrel
{P}{%
\rightarrow}0.
\]
\item[(D6)] For $A \in\{C,L\}$, it holds that
\begin{eqnarray*}
\sup_{y\in\lbrack0,T]} \vert\widehat{S}%
_{x,A}(y)-S_{x}(y) \vert&\stackrel{P}{\rightarrow}&0, \\
 {\sup_{y\in\lbrack0,T]} } \vert{S}^*
_{x,A}(y)-S_{x}(y) \vert&\stackrel{P}{\rightarrow}&0.
\end{eqnarray*}
Here, ${S}^* _{x,A}(y)$ is defined as $\widehat{S}_{x,A}(y)$ in (\ref
{def1S}%
), (\ref{def2S}), (\ref{def1SL}) and (\ref{def2SL}), but with
$\overline
{N}%
_{i}(y)$ replaced by $\overline{\Lambda}_{i}(y)$. (An explicit
definition of
${S}^* _{x,C}(y)$ is also given in the proof of Theorem 3.1.)
\end{enumerate}

These assumptions are rather standard smoothing assumptions. Assumptions
(D4)--(D6) are low-level assumptions. We chose them instead of
high-level assumptions to avoid more specific assumptions on the
censoring. For
the unfiltered case, these assumptions are classical smoothing results. For
the filtered case, consider first the case of random right censoring. Then
\[
n^{-1} \sum_{i=1}^n K_b(x-X_i) C_i(y) Z_i(y) = n^{-1} \sum_{i=1}^n
K_b(x-X_i) I(Y_i^* > y),
\]
where $Y_i^*$ is the minimum of the survival time $Y_i$ and the censoring
time $C_i$, which are supposed to be independent of each other given $X_i$.
It is easily seen that the latter quantity converges to $\varphi_x(y) :=
f(x) P(Y^* > y | X=x)$ uniformly in $x\in R_X$ and $y\in\lbrack0,T]$.
Other examples of filtering (including, e.g., left and/or right truncation
and/or censoring) can be handled in a similar way. Assumption (D5) is only
needed for the asymptotic result based on local constant smoothing and not
for local linear smoothing.

\begin{theorem}
\label{gA} Suppose that assumptions \textup{(D1)--(D6)} hold. There then exist
bounded continuous functions $\beta_{A}$ and $v_{A},$ $A\in\{
C,L \} $, such that for all $x \in R_X$,
\[
\sqrt{nb}\bigl(\hat{g}_{A,\mathrm{mn}}^T(x)-g_{\mathrm{mn}}^T(x)-b^{2}\beta
_{A}(x)\bigr) \Longrightarrow  N(0,v_{A}(x)),
\]
where%
\begin{eqnarray*}
\beta_{C}(x)
&=&\frac{1}{2}\mu_{2}(K)\int_0^T S_{x}(y)\int_{0}^{y}
\biggl\{
\frac{\partial^{2}\alpha_{x}(u)}{\partial x^{2}}+2\frac{\partial
\alpha
_{x}(u)}{\partial x}\frac{\partial\varphi_{x}(u)}{\partial x} \biggr\}
\,\mathrm{d}u\,\mathrm{d}y ,
\\
\beta_{L}(x) &=&\frac{1}{2}\mu_{2}(K)\int_0^T S_{x}(y)\int
_{0}^{y}\frac{%
\partial^{2}\alpha_{x}(u)}{\partial x^{2}}\,\mathrm{d}u\,\mathrm{d}y , \\
v_{C}(x) &=&\|K\|_{2}^{2}\int\frac{\alpha_{x}(u)}{\varphi_{x}(u)}\biggl \{
\int_{u}^T S_{x}(y)\,\mathrm{d}y \biggr\} ^{2}\,\mathrm{d}u , \\
v_{L}(x) &=&v_{C}(x) .
\end{eqnarray*}
\end{theorem}

%In the Appendix we will show that $\widehat{g}_{A,mn}(x)$ reduces to
%the
%local constant and local linear estimators of $g_{mn}(x)$ when the
%data $%
To be consistent with the theory for kernel regression estimators, it must
be that in the absence of filtering,%
\[
v_{C}(x)=\|K\|_{2}^{2}\frac{\sigma^{2}(x)}{f(x)},
\]
where $\sigma^{2}(x)=\operatorname{var}[Y|X = x]$ and $f(x)$ is the
covariate density$.$ Note that
\[
\operatorname{var}[Y|X = x]=2\int uS_{x}(u)\,\mathrm{d}u- \biggl( \int
S_{x}(u)\,\mathrm{d}u \biggr) ^{2}.
\]
In the absence of filtering, $\varphi_{x}(u)=f(x)S_{x}(u).$ Therefore, it
should be the case that%
\[
\int\frac{\alpha_{x}(u)}{S_{x}(u)} \biggl\{ \int_{u}S_{x}(y)\,\mathrm{d}y \biggr\}
^{2}\,\mathrm{d}u=2\int uS_{x}(u)\,\mathrm{d}u- \biggl( \int S_{x}(u)\,\mathrm{d}u \biggr) ^{2}.
\]
This follows by integration by parts. %We just show for a special case.
%Suppose that $Y$ is independent of $X$ with an exponential
%distribution, $%
%f(u)=\lambda\exp(-\lambda u).$ Then we have $\alpha(u)=\lambda$
%and $%
%S(u)=\exp(-\lambda u)$ so that $\int_{u}^{\infty}S(y)dy=\exp(-
%u)/\lambda.$ It follows that
%}S(y)dy \} ^{2}du=\int_{0}^{\infty}\frac{1}{\lambda}\exp(-\lambda
%u)du=\frac{1}{\lambda^{2}}.
%On the other hand, the variance of an exponential random variable is $%
%1/\lambda^{2}$ as required.

For $g_{\mathrm{med}}(x)$, it has been shown in   \cite{Keilegom1998}
that $\hat g_{C,\mathrm{med}}(x)$ is asymptotically normal when the data are
subject to random right censoring. It can be shown that this result
continues to hold true for general filtering patterns.

%s4 ###
\section{Estimation under common shape of the error distribution}\label{sec4}

Under some circumstances, it may be plausible to assume that the error
distribution, when adjusted for the mean or the median, is generated by the
same underlying shape. If there is a substantial level of filtering,
then one
can envision areas where truncation or censoring imply that we do not have
local information on the entire shape of the error distribution around every~$x.$ One
can alleviate this by imposing assumptions on the shape of these
local error distributions. The simplest assumption in this connection is
simply that $\alpha_{\varepsilon_{0}|x}$ does not depend on $x$,
where $%
\varepsilon_{0}$ is the error term in model (\ref{eq:genmodel2}).
This is
%
%e11 ###
\begin{equation}
\alpha_{\varepsilon_{0}|x}(u)\equiv\alpha_{0}(u) \label{mod2}
\end{equation}
for some $\alpha_{0}$ and all $u\geq0.$ If this assumption is true,
then it can be used to improve estimation, even in the case without
filtering, as we now discuss. The notion of efficiency is here tied to
asymptotic variance, which yields mean-squared error holding bias constant,
and comes from the classical parametric theory of likelihood. The local
likelihood method was introduced in
  \cite{Tibshirani1984} and has been applied
in many other contexts. Tibshirani \cite{Tibshirani1984}, Chapter 5, presents the
justification for the local likelihood method (in the context of an
exponential family): the author shows that its asymptotic variance is
the same as
the asymptotic variance of the maximum likelihood estimator (MLE) of a
correctly specified parametric model
at the point of interest using the same number of observations as the local
likelihood method. This type of result has been shown in other
settings, for
example, Linton and Xiao \cite{Linton2002} establish efficiency of a local likelihood
estimator in the context of nonparametric regression with additive errors.
In generalized additive models, Linton \cite{Linton1997,Linton2000} shows the improvement
according to variance obtainable by the local likelihood method.

In what follows, $g(x)$ is either $g_{\mathrm{mn}}(x)$ or $g_{\mathrm{med}}(x)$ and similarly
for the estimators of $g(x)$.

%s4.1 ###
\subsection{Oracle estimation of the location $g(x)$}\label{sec4.1}

First, we note that both the local constant and the local linear kernel
estimator of the full marker-dependent hazard model have the form
\[
\hat{\alpha}_{x,A}(y)=\frac{O_{x,y}^{A}}{E_{x,y}^{A}},
\]
where $A$ equals $C$ for the local constant case and $A$ equals $L$ for the
local linear case. Let us suppose that an oracle told us what $\alpha_{0}$
is. We define the local constant estimator and the local linear estimators
of $g$ based on the assumption (\ref{eq:alfanul}) to be any minimizer $
\hat{g}_{A}^{o}$ of the criterion function
\[
\int\!\!\int\biggl[ \hat{\alpha}_{x,A}(y)-\frac{1}{g(x)}\alpha_{0} \biggl\{
\frac{y}{g(x)} \biggr\} \biggr] ^{2} \{ \hat{\alpha}%
_{x,A}(y) \} ^{-1}E_{x,y}^{A}w(x,y)\,\mathrm{d}x\,\mathrm{d}y,
\]
where $w(x,y)$ is an appropriate weight function. This is motivated by the
theory of minimum chi-squared estimation \cite{Berkson1980}, in which
efficiency is achieved by weighting a least-squares criterion with the
inverse of the asymptotic variance of the unrestricted estimator (in this
case, $\hat{\alpha}_{x,A}(y),$ which has asymptotic variance
$\alpha
_{x}(y)/\varphi_{x}(y)$, where $\varphi_{x}(y)$ is the probability limit
of the exposure $E_{x,y}^{A}).$ For a fixed $x$, this expression is
minimized by minimizing the pointwise criterion%
%
%e12 ###
\begin{equation}
\hat{l}_{\alpha_{0}}(\theta;x)=\int\biggl[ \hat{\alpha}_{x,A}(y)-
\frac{1}{\theta}\alpha_{0} \biggl\{ \frac{y}{\theta} \biggr\} \biggr]
^{2} \{ \hat{\alpha}_{x,A}(y) \} ^{-1}E_{x,y}^{A}w(x,y)\,\mathrm{d}y
\label{t2}
\end{equation}
with respect to $\theta$ and setting $\hat{g}_{A}^{o}(x)=\widehat
{%
\theta}=\arg\min_{\theta\in\Theta}\hat{l}_{\alpha
_{0}}(\theta;x)$
for some compact set $\Theta$ not containing 0. This is a nonlinear
estimator, not obtainable in closed form.

Define%
\[
l(\theta;x)=\int_{0}^{\infty} \biggl[ \alpha_x (y)-\frac{1}{\theta
}\alpha
_{0} \biggl\{ \frac{y}{\theta} \biggr\} \biggr] ^{2} \{ \alpha_x
(y) \} ^{-1}\varphi_{x}(y) w(x,y) \,\mathrm{d}y
\]
and let $\theta_{0}=g(x).$ %\newline

For the asymptotic result below, we need to assume the following:

\begin{enumerate}[(A5)]
\item[(A1)] (i) The weight function $w(x,y)$ is continuous and
satisfies $%
w(x,y)=0$ for $(x,y)\notin I$ and $0\leq w(x,y)\leq a$ for all
$(x,y)\in I$%
, where $0<a<\infty$ and $I=\{(x,y)\dvtx x\in R_{X},\tau_{x}\leq y\leq
T_{x}\}$%
, where $\tau_{x}$ and $T_{x}$ are continuous functions and where, as in
(D1)--(D5), $R_X$ is a bounded interval in the interior of the support
of $X$.

\noindent
(ii) There exists a continuous function $\varphi_{x}(\cdot)$ with $%
\inf_{(x,y)\in I}\varphi_{x}(y)>0$ such that the convergence
statements in
(D4) and (D5) hold with the supremum running over $(x,y)\in I$ instead
of $%
y\in[0,T]$. The function $\varphi_{x}(y)$ is twice continuously
differentiable in $y$ for $(x,y)\in I$.
\item[(A2)] The function $\alpha_{x}(y)=g(x)^{-1}\alpha_{0}[g(x)^{-1}y]$
is twice continuously differentiable in $(x,y)\in I$ and $\inf
_{(x,y)\in
I}\alpha_{x}(y)>0$.
\item[(A3)] The probability density functions $K$ and $k$ are symmetric
around 0 and have support $[-1,1]$, $\int uK(u)\,\mathrm{d}u=\int uk(u)\,\mathrm{d}u=0$,
$\int
u^{2}K(u)\,\mathrm{d}u\neq0$, $\int u^{2}k(u)\,\mathrm{d}u\neq0$, and $K$ and $k$ are twice
continuously differentiable.
\item[(A4)] For all $\varepsilon>0$, $\inf_{|\theta-\theta
_{0}|>\varepsilon}|l(\theta;x)-l(\theta_{0};x)|>0$, $l(\theta,x)$ is
twice differentiable with respect to $\theta$ in a neighborhood of
$\theta_0$
and $l^{\prime\prime}(\theta_{0};x)>0$.
\item[(A5)] The bandwidths $h$ and $b$ satisfy $h\rightarrow0$, $%
b\rightarrow0$, $nhb\rightarrow\infty$, $nh^{4}b=\mathrm{O}(1)$ and $nb^{5}=\mathrm{O}(1)$.
\end{enumerate}

%Then assume that for some $\delta>0$
%l(\theta_{0};x)<\inf_{|\theta-\theta_{0}|>\delta}l(\theta;x).

Conditions (A2), (A3) and (A5) are standard smoothing assumptions.
Assumption (A1) is stated uniformly in $x$ because such a uniform
version is
required in the later Theorems~\ref{alpha0} and~\ref{g2step}.

\begin{theorem}
\label{goracle} Suppose that assumptions \textup{(A1)--(A5)} hold. There then exist
bounded continuous functions $\beta_{A1}^{o}$ and $\beta_{A2}^{o}$,
$A\in
\{ C,L \} $, such that for all $x\in R_{X}$,
\[
\sqrt{nb}\bigl(\hat{g}_{A}^{o}(x)-g(x)-h^{2}\beta
_{A1}^{o}(x)-b^{2}\beta
_{A2}^{o}(x)\bigr)  \Longrightarrow   N(0,v_{A}^{o}(x)),
\]
where, with $s_{0}(u)=1+u\alpha_{0}^{\prime}(u)/\alpha_{0}(u),$%
\begin{eqnarray*}
v_{C}^{o}(x) &=&g(x)^{3}\|K\|_{2}^{2} \biggl[ \int s_{0}^{2} \biggl(\frac{y}{g(x)}%
\biggr)\alpha_{0}\biggl (\frac{y}{g(x)} \biggr)\varphi_{x}(y)w(x,y)\,\mathrm{d}y \biggr] ^{-1} ,
\\
v_{L}^{o}(x) &=&v_{C}^{o}(x) .
\end{eqnarray*}
\end{theorem}

In the absence of filtering, the optimal estimator of $g(x)$, given the
knowledge of $\alpha_{0}(\cdot)$ or, equivalently, of the density
$f_{\varepsilon
}(\cdot)$ of $\varepsilon_{0}$, is the local likelihood estimator that
maximizes%
\[
l(\theta)=\sum_{i=1}^{n}K_{b} ( x-X_{i} )\biggl \{ \ln
f_{\varepsilon} \biggl( \frac{Y_{i}}{\theta} \biggr) -\ln\theta\biggr\} ,
\]
which has score function%
\[
s_{\theta}=-\frac{1}{\theta}\sum_{i=1}^{n}K_{b} ( x-X_{i} )
\biggl\{ \varepsilon_{i}(\theta)\frac{f_{\varepsilon}^{\prime}}{%
f_{\varepsilon}} ( \varepsilon_{i}(\theta) ) +1 \biggr\} ,
\]
where $\varepsilon_{i}(\theta)=Y_{i}/\theta.$ The object
$s_{\varepsilon
}(u)=u(f_{\varepsilon}^{\prime}/f_{\varepsilon})(u)+1$ is known as the
\textit{Fisher scale score} and $I_{2}(f_{\varepsilon})=\int
s_{\varepsilon
}^{2}(u)f_{\varepsilon}(u)\,\mathrm{d}u$ is the corresponding information. One can
show that the asymptotic variance of this oracle local likelihood estimator
is%
%
%e13 ###
\begin{equation}
\frac{\|K\|_{2}^{2}g(x)^{2}}{f(x)I_{2}(f_{\varepsilon})}. \label{jj}
\end{equation}
Supposing that we had $k_{n}= \lfloor nbf(x)/\|K\|_{2}^{2} \rfloor$
observations from the model $Y=g(x)\varepsilon,$ the MLE of $\theta
_{0}=g(x)$ would have asymptotic variance
$g(x)^{2}/I_{2}(f_{\varepsilon
})k_{n}.$ In this sense, the local likelihood method has the efficiency of
the MLE from a sample of size $k_{n}.$

By Efron and Johnstone \cite{Efron1990}, we have%
\[
I_{2}(f_{\varepsilon})=\int\biggl( 1+u\frac{\alpha_{0}^{\prime}(u)}{%
\alpha_{0}(u)} \biggr) ^{2}f_{\varepsilon}(u)\,\mathrm{d}u,
\]
which explains the form of the asymptotic variance above. Suppose that
we take $%
w(x,y)=1$, make a change of variables $y\rightarrow u=y/g(x)$ in $%
v_{C}^{o}(x)$ and make use of the fact that, under no filtering,
$\varphi
_{x}(y)=f(x)S_{x}(y)$ and so $\alpha_{0}(u)\varphi
_{x}(ug(x))=f_{\varepsilon}(u)f(x).$ Then $v_{C}^{o}(x) ={}$(\ref{jj}). This
shows that $\hat{g}_{A}^{o}(x)$ is asymptotically equivalent to the
oracle local likelihood method, that is, efficient in this sense.

%s4.2 ###
\subsection{Estimation with unknown $\protect\alpha_{0}$}\label{sec4.2}

For a given $g(\cdot)$, an estimator of $\alpha_{0}$ can be based on the
minimization principle
\[
\hat{\alpha}_{0,A}^{o}=\operatorname{\arg\min}\limits_{\alpha
(\cdot)}\int\!\!\int\biggl[
\hat{\alpha}_{x,A}(y)-\frac{1}{g(x)}\alpha\biggl\{ \frac{y}{g(x)}%
\biggr\} \biggr] ^{2} \{ \hat{\alpha}_{x,A}(y) \}
^{-1}E_{x,y}^{A} w(x,y) \,\mathrm{d}x\,\mathrm{d}y,
\]
where the choice of weighting function is again motivated by efficiency
considerations. Changing variables $y\mapsto u=y/g(x),$ the objective
function becomes%
\begin{eqnarray*}
&&\int\!\!\int\biggl[ \hat{\alpha}_{x,A}(ug(x))-\frac{1}{g(x)}\alpha(u)%
\biggr] ^{2}g(x) \{ \hat{\alpha}_{x,A}(ug(x)) \}
^{-1}E_{x,ug(x)}^{A} w(x,ug(x)) \,\mathrm{d}x\,\mathrm{d}u \\
&&\quad =\int\!\!\int[ g(x)\hat{\alpha}_{x,A}(ug(x))-\alpha(u) ] ^{2}%
\frac{E_{x,ug(x)}^{A}}{g(x)\hat{\alpha}_{x,A}(ug(x))} w(x,ug(x))
\,\mathrm{d}x\,\mathrm{d}u ,
\end{eqnarray*}
ignoring support considerations. Then, because $\alpha$ does not
depend on $%
x$, we can replace it by the pointwise criteria
\[
\hat{l}_{g}(\alpha;u)=\int[ g(x)\hat{\alpha}%
_{x,A}(ug(x))-\alpha] ^{2}\frac{E_{x,ug(x)}^{A}}{g(x)\hat{\alpha}
_{x,A}(ug(x))} w(x,ug(x)) \,\mathrm{d}x
\]
for each $u,$ whence we obtain the closed form solution%
\[
\hat{\alpha}_{0,A}^{o}(y)=
\frac{
\int E_{x,y g(x)}^{A} w(x,yg(x))
\,\mathrm{d}x
}{%
\int
({E_{x,y g(x)}^{A} w(x,yg(x))}/({g(x)\hat{\alpha}_{x,A}(yg(x))}))
\,\mathrm{d}x
}
.
\]
In practice, one computes $\hat{\alpha}_{0,A}(y)$ as (\ref{eq:alfaet})
with $g(x)$ replaced by a preliminary completely nonparametric
estimator $%
\tilde g$, that is,
%
%e14 ###
\begin{equation}
\hat{\alpha}_{0,A}(y)=\frac{
\int E_{x,y\tilde{g}(x)}^{A}
w(x,y\tilde
g(x)) \,\mathrm{d}x}{\int({E_{x,y \tilde{g}(x)}^{A}w(x,y \tilde
g(x))}/({\tilde
{g}(x)%
\hat{\alpha}_{x,A}(y\tilde{g}(x))}) )\,\mathrm{d}x}. \label{eq:alfaet}
\end{equation}

Let $y$ be a fixed value, that is, such that $\tau\le y \le T$, where\vspace*{-1pt}
$\tau
> \inf_{x \in R_X} \frac{\tau_x}{g(x)}$ and $T < \sup_{x \in R_X}
\frac
{T_x}{%
g(x)}$ (and where we assume that $\inf_{x\in R_X} g(x)>0$). We require the
following assumptions:

\begin{enumerate}[(B3)]
\item[(B1)] The preliminary estimator $\tilde g(\cdot)$ satisfies
$\sup_{x
\in R_X} |\tilde g(x)-g(x)| = \mathrm{O}_P((nb)^{-1/2}\times\break  (\log n)^{1/2})$.

\item[(B2)] The function $g(x)$ is twice continuously differentiable
in $x
\in R_X$ and $\inf_{x \in R_X} g(x) > 0$.

\item[(B3)] The bandwidths $h$ and $b$ satisfy $h \rightarrow0$, $b
\rightarrow0$, $nhb \rightarrow\infty$, $nb^4h = \mathrm{O}(1)$, $nh^5=\mathrm{O}(1)$
and $%
nh^2b(\log n)^{-1} \rightarrow\infty$.
\end{enumerate}

\begin{theorem}
\label{alpha0} Suppose that assumptions \textup{(A1)--(A3)} and \textup{(B1)--(B3)} hold.
There then exist bounded continuous functions $b_{A1}$ and $b_{A2}$,
$A\in
\{ C,L \}$, such that for all $\tau\le y \le T$,
\[
\sqrt{nh}\bigl(\hat{\alpha}_{0,A}(y)-%
\alpha_{0}(y)-h^{2}b_{A1}(y)-b^{2}b_{A2}(y)\bigr)  \Longrightarrow   N(0,s_{A}(y)),
\]
where
\begin{eqnarray*}
 s_{C}(y) &=& \|k\|_{2}^{2} \frac{1}{B^o(y)^2} E \{
E[C(yg(X))|Y=yg(X),X] f_X(yg(X)) w^2(X,yg(X)) \} , \\
  s_{L}(y) &=& s_{C}(y) , \\
  B^o(y) &=& E[(CZ)\{yg(X)\} w\{X,yg(X)\}] / \alpha_0(y) .
\end{eqnarray*}
%
%and
%H_y(x) = \frac{E[(CZ)\{yg(x)\}|X=x] w\{x,yg(x)\}}{E[C\{yg(x)\}|X=x]}.
\end{theorem}

%Note that when the data are subject to random right censoring, then
%$C_i(y)
%= I(\min(Y_i,U_i) \ge y)$ and $(C_iZ_i)(\cdot) \equiv C_i(\cdot)$,
%where $%
%U_i $ is the censoring time, and hence the variance $s_A(y)$ reduces
%in that
%case to
%s_A(y) = \|k\|_2^2 \alpha_0^2(y)\frac{E ([1-G_X\{yg(X)\}] f_X\{yg(X)\}
%w\{X,yg(X)\} )}{E^2 ([1-H_X\{yg(X)\}] w\{X,yg(X)\} )},
%where $G_x(y) = P(U \le y|X=x)$ and $H_x(y) = P(\min(Y,U) \le y|X=x)$.

Finally, we compute a new estimate of $g$ using the estimate of $\alpha
_{0}. $ Specifically, define the weighted least-squares objective function
\[
\hat{l}_{\hat{\alpha}}(\theta;x )=\int\biggl[ \hat{\alpha}%
_{x,A}(y)-\frac{1}{\theta}\hat{\alpha}_{0} \biggl\{ \frac{y}{\theta}
\biggr\} \biggr] ^{2} \{ \hat{\alpha}_{x,A}(y) \}
^{-1}E_{x,y}^{A} w(x,y) \,\mathrm{d}y
\]
with $A=C$ or $A=L$, and with $\hat{\alpha}_{0}$ equal to
$\hat
{%
\alpha}_{0,C}$, $\hat{\alpha}_{0,L}$ or another estimator of
${\alpha}%
_{0}$. Then let
\[
\hat{g}^{2-\mathit{step}}(x)=\arg\min_{\theta\in I_n(x) }\hat
{l}_{\hat{%
\alpha}}(\theta;x ),
\]
where the argmin runs over a shrinking neighborhood $I_n(x)$ of a consistent
estimator of $g(x)$. In the next theorem, we state that under some conditions
on the estimator $\hat{\alpha}_{0}$, we obtain the same variance and
bias as in the oracle case. One possibility is to use the estimator of $g$
given in Section~\ref{sec3} as preliminary estimator and to base the final estimation
of $g$ on the method of the above Section~\ref{sec4.1}, but replacing the oracle
$\alpha_{0}$ by $\hat\alpha_{0,A}$. We make use of the following
additional assumptions:

\begin{enumerate}[(C3)]
\item[(C1)] For a neighborhood $J(x)$ of the closed interval $[\tau
_{x}/g(x),T_{x}/g(x)]$, it holds uniformly for $z\in J(x)$ that
\begin{eqnarray*}
 \hat{\alpha}_{0}(z)-{\alpha}_{0}(z)&=&\mathrm{O}_{P}(\delta_{0,n}), \\
\hat{\alpha}_{0}^{\prime}(z)-{\alpha}_{0}^{\prime
}(z)&=&\mathrm{O}_{P}(\delta
_{1,n}), \\
 \hat{\alpha}_{0}^{\prime\prime}(z)-{\alpha}_{0}^{\prime
\prime
}(z)&=&\mathrm{O}_{P}(\delta_{2,n})
\end{eqnarray*}
for sequences $\delta_{0,n}$, $\delta_{1,n}$ and $\delta_{2,n}$ with
$%
\delta_{0,n}=\mathrm{o}((\log n)^{-1/2}h^{1/2})$, $\delta_{1,n}=\mathrm{o}(1)$,
$\delta_{2,n}=\mathrm{o}((nbh)^{1/2}(\log n)^{-1/2})$, $\delta_{0,n}\delta
_{2,n}=\mathrm{o}(1)$, $\delta_{1,n}\delta_{2,n}=\mathrm{o}(1)$ and $\delta
_{0,n}\delta
_{1,n}=\mathrm{o}(n^{-2/5})$.
\item[(C2)] With a bounded function $\gamma(x)$, it holds that
\[
\int[ \hat{\alpha}(y)-\alpha_{0}(y) ] \rho_{0}\biggl \{
\frac{y}{g(x)} \biggr\} \frac{1}{g(x)^{2}}\frac{E_{x,y}^{C}}{\alpha
_{x}(y)}%
w(x,y)\,\mathrm{d}y-h^{2}\gamma(x)=\mathrm{o}_{P}\bigl((nb)^{-1/2}\bigr),
\]
where $\rho_{0}(u)=\alpha_{0}(u)+u\alpha_{0}^{\prime}(u)$.

\item[(C3)] The bandwidths $h$ and $b$ satisfy $nhb\rightarrow\infty
$, $%
nh^{5}=\mathrm{O}(1)$, $nb^{5}=\mathrm{O}(1)$ and $1/(nb^{5})=\mathrm{O}(1)$.
\end{enumerate}

These assumptions are rather weak. Assumption (C1) is fulfilled for a
standard one-dimensional kernel smoother which fulfills the conditions
with $%
\delta_{0,n}=(\log n)^{1/2}n^{-2/5}$, $\delta_{1,n}=(\log n)^{1/2}n^{-1/5}$
and $\delta_{2,n}=(\log n)^{1/2}$. The assumption is fulfilled under much
slower rates of convergence. The assumption could be replaced by another
type of condition using the general approach of Mammen and Nielsen \cite{Mammen2007}
based on cross-validation arguments. Assumption (C2) is a standard property
of kernel smoothers: kernel smoothers are local weighted averages.
Integration of the estimator leads to a global weighted average with
stochastic part of parametric rate $n^{-1/2}$. Typically, the rate of the
bias part does not change.

\begin{theorem}
\label{g2step} Suppose that assumptions \textup{(A1)--(A5)} and \textup{(C1)--(C3)} hold.
There then exist bounded continuous functions $\beta^{2-\mathit{step}}$ and $v^{2-\mathit{step}}$
such that for all $x \in R_X$,
\[
\sqrt{nb}\bigl(\hat{g}^{2-\mathit{step}}(x)-g(x)-b^{2}\beta
^{2-\mathit{step}}(x)\bigr)  \Longrightarrow   N(0,v^{2-\mathit{step}}(x)),
\]
where $v^{2-\mathit{step}}(x)=v_{C}^{o}(x)=v_{L}^{o}(x).$
\end{theorem}

This shows that the two-step estimator achieves the desired oracle property.

%s5 ###
\section{Numerical results}\label{sec5}

In this section, we look at the small-sample performance of our estimators.
The design involves a combination of commonly occurring features in the
literature:
we take the true underlying regression function to be identical to that of
Fan and Gijbels \cite{Fan1994}, but our disturbance term has a different
distribution and we also consider a different censoring mechanism. Thus,
\begin{eqnarray*}
Y_{i} &=&g_{\mathrm{mn}} ( X_{i} ) \varepsilon_{i} , \\
g_{\mathrm{mn}} ( x ) &=&4.5-64x^{2} ( 1-x ) ^{2}-16 (
x-0.5 ) ^{2},
\end{eqnarray*}
where $X_{i}\sim U [ 0,1 ] ,$ $\varepsilon_{i}\sim U [ 0.5,1.5%
] ,$ while $X_{i}$ and $\varepsilon_{i}$ are independent and $%
E(\varepsilon_i)=1$. The censoring time mechanism is independent of the
covariate and constructed as follows:%
\[
U_{i}= \cases{
V_{i} &\quad  if $W_{i}<0.5$, \cr
+\infty&\quad otherwise,
}
\]
where $V_{i}\sim \mathit{Beta} ( 1,3 ) ,$ $W_{i}\sim \mathit{Beta} (
1,0.75 ) $ and we observe $\{Y_{i}\wedge U_{i},\delta_{i}=1 (
Y_i<U_i ) ,X_i\}$, that is, an example of right censoring.\vadjust{\goodbreak}

We employ two methods of estimation of $g_{\mathrm{mn}} ( X ) $: the simple
local constant estimation of Section~\ref{sec3} and the feasible oracle
estimation, as
discussed in Section~\ref{sec4.2}. For the purposes of illustration, we use Silverman's
rule of thumb bandwidth and the built-in minimization routine based on the
golden section search and parabolic interpolation. For the more efficient
estimator, we note that using the one-dimensional grid search gives a very
similar estimate.\looseness=1

%f1 ###
\begin{figure}

\includegraphics{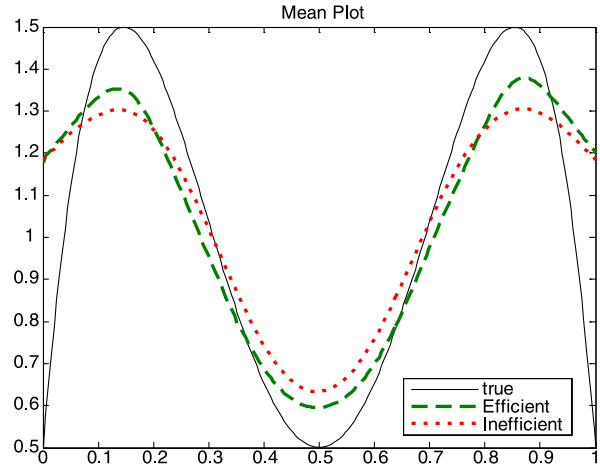}

  \caption{Plot of the mean of
the estimated regression curve.}\label{fig1}
\end{figure}
%f2 ###
\begin{figure}

\includegraphics{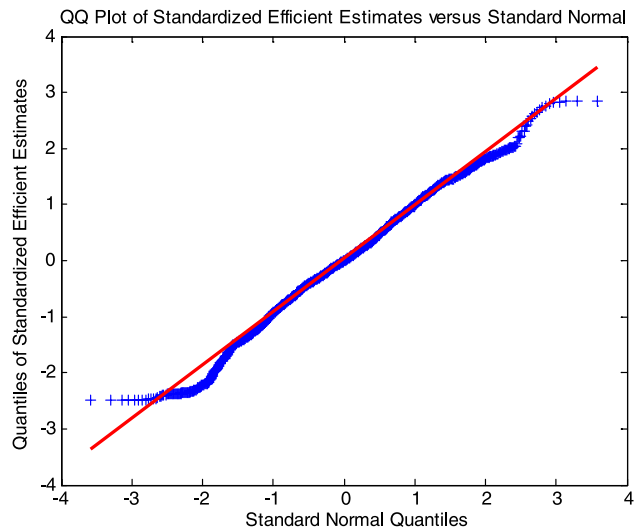}

  \caption{QQ-plot of
standardized efficient estimates versus standard normal.}\label{fig2}
\end{figure}

\begin{figure}[b]

\includegraphics{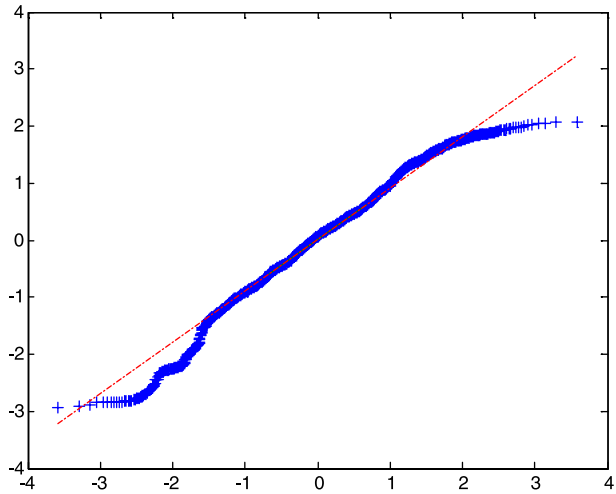}

  \caption{QQ-plot of
standardized inefficient estimates versus standard normal.}\label{fig3}
\end{figure}

\begin{figure}

\includegraphics{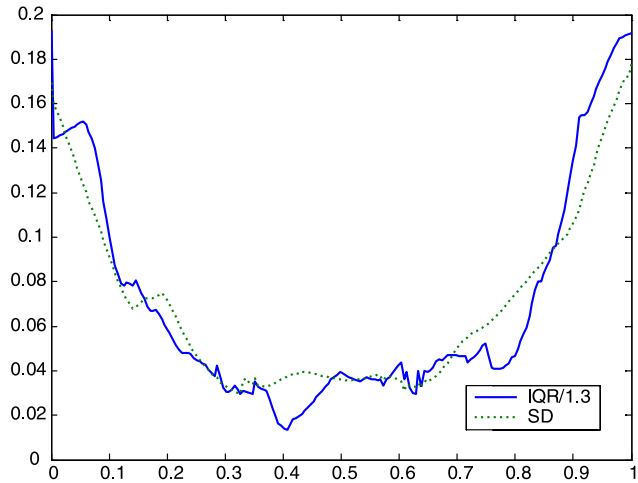}

  \caption{Normalized
interquartile range and standard deviation plots for efficient estimates.}\label{fig4}\vspace*{-3pt}
\end{figure}

We use a sample of size $250$ and $15$ replications over $200$ evenly
spaced grids on $ [ 0,1 ] $. In this example, approximately $25\%$
of the $200$ observations are censored. Figure~\ref{fig1} displays the average (over
replications) of the two estimates. The true regression function chosen
possesses a high degree of curvature, with the function increasing less
steep to the right of $0.5$ than to the left of $0.5$. Both estimates are
capable of capturing the basic structure of the true curve. The efficient
estimate appears to adapt better at both peaks and troughs, and the quality
of fit declines with the steepness of the true curve. Although it is
not shown
here, the relative performance of the simple local constant estimator
improves toward the feasible oracle estimates when the true regression
function has lower degree variation. Figures~\ref{fig2} and~\ref{fig3} are the QQ-plots  for
the  efficient and inefficient estimates, respectively (i.e., $(\hat
{g}%
-E\hat{g})/\mathit{std}(\hat{g}))$. The linear trends in   the QQ-plots are
distinct with the efficient estimates performing a little better away from
the sample means. Figure~\ref{fig4} plots the interquartile range (divided by 1.3)
and the standard deviation (across replications) for the efficient estimate
against grid points. Performance clearly worsens in the boundary region.

%f3 ###

%f4 ###

Since it is widely perceived that the Silverman's rule of thumb bandwidth
tends to oversmooth, we also performed some experiments with smaller bandwidths.
Smaller bandwidth leads to much larger simulation time during optimization,
due to higher variance. In terms of goodness of fit, it does not make a
big difference with the
feasible oracle estimation. However, the improvement of fit for the simple
local constant estimation is more pronounced. In that case, the feasible
oracle estimation still performs better than the simple estimator, as
expected.

%sA ###
\begin{appendix}\label{appm}

%sB ###
\section*{Appendix}

%{section}{1} \setcounter{definition}{0} \setcounter{subsection}{0}

\renewcommand{\theequation}{\arabic{equation}}

%sB.1 ###
\subsection{Local linear estimation}

In this section, we first define the local linear marker-dependent estimator,
$\hat{\alpha}_{x,L}(y),$ as defined in   \cite{Nielsen1998}, page 118,
%
%eB.1 ###
\begin{equation}
\hat{\alpha}_{x,L}(y)=\frac{O_{x,y}^{L}}{E_{x,y}^{L}}, \label{def1SL}
\end{equation}
where, with $w=(w_j)_{j=1}^{d+1}=(x,y)$ and
$W_{i}(u)=(W_{ij}(u))_{j=1}^{d+1}=(X_{i},u)$ (to simplify the notation, we
consider the same kernel and bandwidth for $x$ and $y$),
\begin{eqnarray} \label{def2SL}
O_{x,y}^{L}&=&n^{-1}\sum_{i=1}^{n}\int{K}_{w,b} \{ w-W_{i}(u) \}\, \mathrm{d}%
\overline{N}_{i}(u) ,\qquad
E_{x,y}^{L}=c_{0}-c_{1}^{T}D^{-1}c_{1} ,
\nonumber\\[-0.5pt]
K_{w,b}(v)&=& \{ K_{b}(v)-K_{b}(v)v^{T}D^{-1}c_{1} \}\qquad
(v \in R^{d+1}),
\nonumber
\\[-0.5pt]
c_{0}&=&n^{-1}\sum_{i=1}^{n}\int K_{b} \{ w-W_{i}(u) \} C_{i}(u)
Z_i(u) \,\mathrm{d}u ,
\\[-0.5pt]
c_{1j}&=&n^{-1}\sum_{i=1}^{n}\int K_{b} \{ w-W_{i}(u) \} \{
w_{j}-W_{ij}(u) \} C_{i}(u) Z_i(u) \,\mathrm{d}u ,
\nonumber\\[-0.5pt]
d_{jk}&=&n^{-1}\sum_{i=1}^{n}\int K_{b} \{ w-W_{i}(u) \} \{
w_{j}-W_{ij}(u) \} \{ w_{k}-W_{ik}(u) \} C_{i}(u) Z_i(u) \,\mathrm{d}u ,
\nonumber
\end{eqnarray}
and $c_1 = (c_{1j})_{j=1}^{d+1}$ and $D = (d_{jk})_{j,k=1}^{d+1}$. We
then consider the local linear estimator of the integrated conditional hazard
function, obtained when we undersmooth in the $y$-direction. First, we
define the necessary kernel constants:
\begin{eqnarray*}
\overline{K}_{x,b}(v)&=& \{ K_{b}(v)-K_{b}(v)v^{T}\overline{D}^{-1}%
\overline{c}_{1} \} /(\overline{c}_{0}-\overline{c}_{1}^{T}\overline
{D}%
^{-1}\overline{c}_{1})\qquad (v \in R^d) ,
\\
\overline{c}_{0}&=&n^{-1}\sum_{i=1}^{n}K_{b} ( x-X_{i} ) ,
\\
\overline{c}_{1j}&=&n^{-1}\sum_{i=1}^{n}K_{b} ( x-X_{i} ) (
x_{j}-X_{ij} ) ,
\\
\overline{d}_{jk}&=&n^{-1}\sum_{i=1}^{n}K_{b} ( x-X_{i} ) (
x_{j}-X_{ij} ) ( x_{k}-X_{ik} ) .
\end{eqnarray*}
We then get that the local linear estimator of the integrated hazard is
\[
\widehat{A }_{x,L}(y)=\int_{0}^{y}\frac{\sum_{i=1}^{n}\overline{K}%
_{x,b} ( x-X_{i} ) }{\sum_{j=1}^{n}\overline{K}_{x,b} (
x-X_{j} ) C_{j}(u)Z_{j}(u)}\,\mathrm{d}\overline{N}_{i}(u) .
\]
We estimate correspondingly the conditional survival function $S_x(y)$ and
the regression functions $g_{\mathrm{mn}}^T(x)$ and $g_{\mathrm{med}}(x)$ by
\begin{eqnarray*}
\widehat{S}_{x,L}(y)&=&{\prod_{0\leq w\leq y\leq T} } \{ 1-d%
\widehat{A }_{x,L}(w) \} ,
\\
\hat{g}_{L,\mathrm{mn}}^T(x)&=&-\int_{0}^{T}y\,\mathrm{d}\widehat{S}_{x,L}(y) ,
\\
\hat{g}_{L,\mathrm{med}}(x)&=&\widehat{S}_{x,L}^{-1}(0.5) .
\end{eqnarray*}

%sB.2 ###
\subsection{Proof of results}

We restrict attention in the proofs to the case of local constant smoothing
(i.e., when $A=C$). The case of local linear smoothing ($A=L$) can be
considered in a very similar way and is therefore omitted. Throughout this
section, we use the notation $A_n \simeq B_n$ to indicate that $A_n = B_n
(1+\mathrm{o}_P(1))$.

First, we state a useful lemma. Its simple proof is omitted. Let
$h(y-)$ be
the limit from the left at $y$ for any cadlag function $h$.

\begin{lemma}
\label{lem} Suppose $A _{1}$ and $A _{2}$ are cadlag functions. Let $%
S_{1}(y)=\prod_{w\leq y}  \{ 1-\mathrm{d}A _{1}(w) \} $, $%
S_{2}(y)=\prod_{w\leq y}  \{ 1-\mathrm{d}A _{2}(w) \} $ and
\[
Q(y)=\frac{S_{1}(y-)}{S_{2}(y-)}-1\mathit{.}
\]
Then\vspace*{-2pt}
\[
\mathrm{d}Q(y)=\frac{S_{1}(y-)}{S_{2}(y-)}\,\mathrm{d} ( A _{1}-A _{2} ) (y).
\]
\end{lemma}

\begin{pf*}{Proof of Theorem \ref{gA}} Define\vspace*{-2pt}
\[
A _{x,C}^{\ast}(y)=\sum_{i=1}^{n}K_{b} ( x-X_{i} ) \int_{0}^{y}%
\frac{\mathrm{d}\overline{\Lambda}_{i}(u)}{%
\sum_{j=1}^{n}K_{b}(x-X_{j})C_{j}(u)Z_{j}(u)}.
\]
Then\vspace*{-2pt}
\[
\widehat{A }_{x,C}(y)-A _{x}^{\ast}(y)=\sum_{i=1}^{n}K_{b} (
x-X_{i} ) \int_{0}^{y}\frac{\mathrm{d}\overline{M}_{i}(u)}{%
\sum_{j=1}^{n}K_{b}(x-X_{j})C_{j}(u)Z_{j}(u)}.
\]

Let $S_{x}^{\ast}(y)={\prod_{w\leq y} } \{ 1-\mathrm{d}A _{x}^{\ast
}(w) \} .$ We then divide our analysis into an analysis of the variable
part\vspace*{-2pt}
%
%eB.3 ###
\begin{equation}
V_{x}(y)=\widehat{S}_{x,C}(y)-S_{x}^{\ast}(y)
\end{equation}
and of the stable part\vspace*{-2pt}
%
%eB.4 ###
\begin{equation}
B_{x}(y)=S_{x}^{\ast}(y)-S_{x}(y).
\end{equation}

Note that $V_{x}(y)=S_{x}^{\ast}(y)Q_{x}^{V}(y),$ where $Q_{x}^{V}(y)=%
\widehat{S}_{x,C}(y)/S_{x}^{\ast}(y)-1,$ $B_{x}(y)=S_{x}(y)Q_{x}^{B}(y)$
and $Q_{x}^{B}(y)=S_{x}^{\ast}(y)/S_{x}(y)-1.$ Using integration by
parts, we obtain\vspace*{-2pt}%\vspace*{-2pt}
%
%_{x,C}(\mathrm{d}y)-S_{x}(\mathrm{d}y)]\\&=&\int_0^T [\widehat
%{S}_{x,C}(y)-S_{x}(y)]\,\mathrm{d}y=\mathcal{V}%
%(x)+\mathcal{B}(x),
\[
\hat{g}^T_{\mathrm{mn}}(x)-g^T_{\mathrm{mn}}(x)=-\int_0^T y[\widehat{S}%
_{x,C}(\mathrm{d}y)-S_{x}(\mathrm{d}y)]=\int_0^T [\widehat
{S}_{x,C}(y)-S_{x}(y)]\,\mathrm{d}y=\mathcal{V}%
(x)+\mathcal{B}(x),
\]
where $\mathcal{V}(x)=\int_0^T V_{x}(y)\,\mathrm{d}y$ and $\mathcal{B}(x)=\int_0^T
B_{x}(y)\,\mathrm{d}y.$ By Lemma \ref{lem}, we have\vspace*{-2pt}%
\begin{eqnarray*}
\mathcal{V}(x) &=&\int_{0}^{T}S_{x}^{\ast}(y)\int
_{0}^{y}\mathrm{d}Q_{x}^{V}(u)\,\mathrm{d}y \\[-1pt]
&=&\int_{0}^{T}S_{x}^{\ast}(y)\int_{0}^{y}\frac{\widehat{S}_{x,C}(u-)}{
S_{x}^{\ast}(u-)}\,\mathrm{d} \{ \widehat{A }_{x,C}(u)-A _{x,C}^{\ast}(u) \}
\,\mathrm{d}y \\[-1pt]
&=&\int_{0}^{T}S_{x}^{\ast}(y)\int_{0}^{y}\frac{\widehat{S}_{x,C}(u-)}{
S_{x}^{\ast}(u-)}\sum_{i=1}^{n}K_{b} ( x-X_{i} )
\frac{\mathrm{d}\overline{M%
}_{i}(u)}{\sum_{j=1}^{n}K_{b}(x-X_{j})C_{j}(u)Z_{j}(u)}\,\mathrm{d}y
\\[-1pt]
%&&\hspace*{118pt}{}\times \frac{\mathrm{d}\overline{M
%}_{i}(u)}{\sum_{j=1}^{n}K_{b}(x-X_{j})C_{j}(u)Z_{j}(u)}\,\mathrm{d}y \\
&=&\sum_{i=1}^{n}\int_{0}^{T}\hat{h}_{x}^{i}(u)\,\mathrm{d}\overline{M}_{i}(u),
\end{eqnarray*}
where\vspace*{-2pt}
\begin{eqnarray*}
\hat{h}_{x}^{i}(u) &=&\int_{u}^{T}S_{x}^{\ast}(y)\frac{\widehat
{S}%
_{x,C}(u-)}{S_{x}^{\ast}(u-)}K_{b} ( x-X_{i} ) \frac{1}{%
\sum_{j=1}^{n}K_{b}(x-X_{j})C_{j}(u)Z_{j}(u)}\,\mathrm{d}y \\
&=&\frac{\widehat{S}_{x,C}(u-)}{S_{x}^{\ast}(u-)}\frac{K_{b} (
x-X_{i} ) }{\sum_{j=1}^{n}K_{b}(x-X_{j})C_{j}(u)Z_{j}(u)}
\int_{u}^{T}S_{x}^{\ast}(y)\,\mathrm{d}y.\
\end{eqnarray*}

Let
\begin{eqnarray*}
\widetilde{\mathcal{V}}(x)
&=&\sum_{i=1}^{n}\int
_{0}^{T}h_{x}^{i}(u)\,\mathrm{d}%
\overline{M}_{i}(u), \\
h_{x}^{i}(u) &=&n^{-1}\frac{K_{b} ( x-X_{i} ) }{\varphi_{x}(u)}%
\int_{u}^{T}S_{x}(y)\,\mathrm{d}y.
\end{eqnarray*}
Then $\mathcal{V}(x)=\widetilde{\mathcal{V}}(x)+\mathrm{o}_{p}(1)$ and by
Nielsen and Linton \cite{Nielsen1995}, Proposition
1, $(nb)^{1/2} \widetilde{\mathcal{V}}%
(x)\Longrightarrow N(0,v(x)),$ where%
\begin{eqnarray*}
v(x) &=&p\lim_{n\rightarrow\infty} nb\sum_{i=1}^{n}\int
h_{x}^{i}(u)^{2}\,\mathrm{d} \langle\overline{M}_{i}(u) \rangle\\
&=&p\lim_{n\rightarrow\infty}m n^{-1} b
\sum_{i=1}^{n}K_{b} ( x-X_{i} ) ^{2}
\int\frac{1}{\varphi
_{x}(u)^{2}} \biggl\{ \int_{u}^{T}S_{x}(y)\,\mathrm{d}y \biggr\} ^{2}\alpha_{X_{i}}
(u)C_{i}(u)Z_{i}(u)\,\mathrm{d}u \\
&=&\|K\|_{2}^{2}\int\frac{\alpha_x (u)}{\varphi_{x}(u)}\biggl \{
\int_{u}^{T}S_{x}(y)\,\mathrm{d}y \biggr\} ^{2}\,\mathrm{d}u.
\end{eqnarray*}
The results given in Theorem \ref{gA} on the variable part follow from
standard martingale theory; see, among many others,   \cite{Nielsen1995}.

We now turn to the bias. Using (D4)--(D6), we have%
\begin{eqnarray*}
\mathcal{B}(x)
&=&\int_{0}^{T}S_{x}(y)\int_{0}^{y}\mathrm{d}Q_{x}^{B}(u)\,\mathrm{d}y \\
&=&\int_{0}^{T}S_{x}(y)\int_{0}^{y}\frac{{S}_{x}^{\ast}(u-)}{S_{x}(u-)}
\,\mathrm{d} \{ A _{x}^{\ast}(u)-A _{x}(u) \} \,\mathrm{d}y \\
&=&\int_{0}^{T}S_{x}(y)\int_{0}^{y}\frac{{S}_{x}^{\ast}(u-)}{S_{x}(u-)}
\sum_{i=1}^{n}K_{b} ( x-X_{i} ) \frac{C_{i}(u)Z_{i}(u) \{
\alpha_{X_{i}} (u)-\alpha_x (u) \} \,\mathrm{d}u}{%
\sum_{j=1}^{n}K_{b}(x-X_{j})C_{j}(u)Z_{j}(u)}\,\mathrm{d}y \\
&=&\frac{1}{2}\mu_{2}(K)b^{2}\int_{0}^{T}S_{x}(y)\int_{0}^{y} \biggl\{
\frac{
\partial^{2}\alpha_x (u)}{\partial x^{2}}+\frac{1} {2} \frac
{\partial
\alpha_x (u)}{\partial x}\frac{\partial\varphi_{x}(u)}{\partial x}
\biggr\}
\,\mathrm{d}u\,\mathrm{d}y+\mathrm{o}_{P}(b^{2}).
\end{eqnarray*}

The derivation of the asymptotic theory of the local linear case parallels
the local constant case. While the variable part has the same asymptotic
distribution, the stable part changes due to the bias properties of the
local linear hazard estimator. By checking the derivation of the stable part
of the local linear kernel hazard estimation of Nielsen \cite{Nielsen1998}, page
119, it
is easy to see that the stable part of the local linear estimator can be
written as
\[
b^{L}(x)=\frac{1}{2}\mu_{2}(K)b^{2}\int_{0}^{T}S_{x}(y)\int
_{0}^{y}\frac{%
\partial^{2}\alpha_x (u)}{\partial x^{2}}\,\mathrm{d}u\,\mathrm{d}y.
\]
\upqed
\end{pf*}

\begin{pf*}{Proof of Theorem \ref{goracle}} Consistency of $\widehat{%
\theta}$ follows from condition (A1) and the fact that%
%
%eB.5 ###
\begin{equation}
\sup_{\theta\in\Theta} \vert\hat{l}_{\alpha_{0}}(\theta
;x)-l(\theta;x) \vert\stackrel{P}{\rightarrow}0 \label{Q1}
\end{equation}
(see, e.g.,  \cite{Vaart1998}, Theorem 5.7, page~45). The result
(\ref{Q1}%
) follows from assumption (A2) and the uniform consistency of $\hat{
\alpha}_{x,C}(y);$ this is established in  \cite{Nielsen1995}, Theorem
2. Actually,
\begin{eqnarray*}
&&\hat{l}_{\alpha_{0}}(\theta;x)-l(\theta;x)\\
&&\quad =\int\biggl[ \hat{%
\alpha}_{x,C}(y)-\frac{1}{\theta}\alpha_{0} \biggl\{ \frac{y}{\theta}%
\biggr\} \biggr] ^{2} \{ \hat{\alpha}_{x,C}(y) \}
^{-1}E_{x,y}^{C} w(x,y) \,\mathrm{d}y
\\
&&\qquad{}-\int\biggl[ \alpha_{x}(y)-\frac{1}{\theta}\alpha_{0} \biggl\{ \frac{y}{%
\theta} \biggr\} \biggr] ^{2} \{ \alpha_{x}(y) \} ^{-1}\varphi
_{x}(y) w(x,y) \,\mathrm{d}y
 \\
&&\quad =\int[ \hat{\alpha}_{x,C}(y)-\alpha_{x}(y) ] ^{2} \{
\alpha_{x}(y) \} ^{-1}\varphi_{x}(y) w(x,y) \,\mathrm{d}y \\
&&\qquad{}+2\int[ \hat{\alpha}_{x,C}(y)-\alpha_{x}(y) ] \biggl[
\alpha_{x}(y)-\frac{1}{\theta}\alpha_{0} \biggl\{ \frac{y}{\theta} \biggr\} %
\biggr] \{ \alpha_{x}(y) \} ^{-1}\varphi_{x}(y) w(x,y) \,\mathrm{d}y \\
&&\qquad{}
+
\int\biggl[ \hat{\alpha}_{x,C}(y)-\frac{1}{\theta}\alpha
_{0} \biggl\{ \frac{y}{\theta} \biggr\} \biggr] ^{2} \\
&&\qquad\hphantom{+\int[}{}\times[ \{ \hat{%
\alpha}_{x,C}(y) \} ^{-1}E_{x,y}^{C}- \{ \alpha_{x}(y) \}
^{-1}\varphi_{x}(y) ] w(x,y) \,\mathrm{d}y
\end{eqnarray*}
and this converges to zero in probability, uniformly in $\theta\in
\Theta$.
%Then use Cauchy-Schwarz inequality etc.

We next establish asymptotic normality. First, we consider the Taylor
expansion%
%
%eB.6 ###
\begin{equation}
0=\hat{l}_{\alpha_{0}}^{\prime}(\widehat{\theta};x)=\hat
{l}%
_{\alpha_{0}}^{\prime}(\theta_{0};x)+\hat{l}_{\alpha
_{0}}^{\prime
\prime}(\theta^{\ast};x)(\widehat{\theta}-\theta_{0}), \label{te}
\end{equation}
where $\theta^{\ast}$ lies between $\widehat{\theta}$ and $\theta_{0}$.
We have
\[
\hat{l}_{\alpha_{0}}^{\prime}(\theta;x)=2\int\biggl[ \hat
{\alpha}
_{x,C}(y)-\frac{1}{\theta}\alpha_{0}\biggl \{ \frac{y}{\theta} \biggr\} %
\biggr] \rho_{0} \biggl\{ \frac{y}{\theta} \biggr\} \frac{1}{\theta^{2}}%
\frac{E_{x,y}^{C}}{\hat{\alpha}_{x,C}(y)} w(x,y) \,\mathrm{d}y,
\]
where $\rho_{0}(u)=\alpha_{0}(u)+u\alpha_{0}^{\prime}(u)$ and
\begin{eqnarray*}
\hat{l}_{\alpha_{0}}^{\prime\prime}(\theta;x )
&=&2\int\rho
_{0}^{2} \biggl\{ \frac{y}{\theta} \biggr\} \frac{1}{\theta^{4}}\frac{%
E_{x,y}^{C}}{\hat{\alpha}_{x,C}(y)} w(x,y) \,\mathrm{d}y
\\
&&{}-2\int\biggl[ \hat{\alpha}_{x,C}(y)-\frac{1}{\theta}\alpha
_{0} \biggl\{ \frac{y}{\theta} \biggr\} \biggr] \rho_{0}^{\prime} \biggl\{
\frac{y}{\theta} \biggr\} \frac{y}{\theta^{4}}\frac
{E_{x,y}^{C}}{\hat{%
\alpha}_{x,C}(y)} w(x,y) \,\mathrm{d}y \\
&&{}-4\int\biggl[ \hat{\alpha}_{x,C}(y)-\frac{1}{\theta}\alpha
_{0} \biggl\{ \frac{y}{\theta} \biggr\} \biggr] \rho_{0} \biggl\{ \frac{y}{%
\theta} \biggr\} \frac{1}{\theta^{3}}\frac{E_{xy}^{C}}{\hat{\alpha}%
_{x,C}(y)} w(x,y) \,\mathrm{d}y.
\end{eqnarray*}

We first establish the properties of $\hat{l}_{\alpha
_{0}}^{\prime
}(\theta_{0};x).$ Recall from   \cite{Nielsen1995} that
%
%eB.7 ###
\begin{equation}
\hat{\alpha}_{x,C}(y)-\alpha_{x}(y)=\frac{\mathcal
{V}_{x,y}+\mathcal{B}%
_{x,y}}{E_{x,y}^{C}}, \label{NL}
\end{equation}
where
\begin{eqnarray*}
\mathcal{V}_{x,y} &=&\frac{1}{n}\sum_{i=1}^{n}\int
K_{b}(x-X_{i})k_{h}(y-y^{\prime})\,\mathrm{d}M_{i}(y^{\prime}) , \\
\mathcal{B}_{x,y} &=&\frac{1}{n}\sum_{i=1}^{n}\int
K_{b}(x-X_{i})k_{h}(y-y^{\prime}) [ \alpha_{X_{i}}(y^{\prime
})-\alpha
_{x}(y) ] C_{i}(y^{\prime}) Z_{i}(y^{\prime})\,\mathrm{d}y^{\prime} .
\end{eqnarray*}
Therefore,
\begin{eqnarray*}
\hat{l}_{\alpha_{0}}^{\prime}(\theta_{0};x) &=&2\int[ \hat{
\alpha}_{x,C}(y)-\alpha_{x}(y) ] \rho_{0} \biggl\{ \frac{y}{g(x)}%
\biggr\} \frac{1}{g(x)^{2}}\frac{E_{x,y}^{C}}{\hat{\alpha}_{x,C}(y)}
w(x,y) \,\mathrm{d}y \\
&=&2\int\rho_{0} \biggl\{ \frac{y}{g(x)} \biggr\} \frac{1}{g(x)^{2}}\frac{%
\mathcal{V}_{x,y}}{\hat{\alpha}_{x,C}(y)} w(x,y) \,\mathrm{d}y\\
&&{}+
2\int\rho
_{0} \biggl\{ \frac{y}{g(x)} \biggr\} \frac{1}{g(x)^{2}}\frac{\mathcal{B}_{x,y}%
}{\hat{\alpha}_{x,C}(y)} w(x,y) \,\mathrm{d}y \\
&\simeq&
2\int\rho_{0} \biggl\{ \frac{y}{g(x)} \biggr\} \frac{1}{g(x)^{2}}%
\frac{\mathcal{V}_{x,y}}{\alpha_{x}(y)} w(x,y) \,\mathrm{d}y\\
&&{}+2\int\rho_{0} \biggl\{
\frac{y}{g(x)} \biggr\} \frac{1}{g(x)^{2}}\frac{\mathcal{B}_{x,y}}{\alpha
_{x}(y)} w(x,y) \,\mathrm{d}y ,
\end{eqnarray*}
where the last line follows from
\cite{Linton2003}, Lemma 3. Consider%
\begin{eqnarray*}
&& \int\rho_{0} \biggl\{ \frac{y}{g(x)} \biggr\} \frac{1}{g(x)^{2}}\frac{%
\mathcal{V}_{x,y}}{\alpha_{x}(y)} w(x,y) \,\mathrm{d}y \\
&&\quad  = \frac{1}{n}\sum_{i=1}^{n}K_{b}(x-X_{i})\int\biggl[ \int\rho
_{0} \biggl\{ \frac{y}{g(x)} \biggr\} \frac{1}{g(x)^{2}}\frac{1}{\alpha_{x}(y)}%
k_{h}(y-y^{\prime}) w(x,y) \,\mathrm{d}y \biggr] \,\mathrm{d}M_{i}(y^{\prime}) \\
&&\quad \simeq\sum_{i=1}^{n}\int h_{ni}(x,u)\,\mathrm{d}M_{i}(u),
\end{eqnarray*}
where $h_{ni}(x,u)=n^{-1}K_{b}(x-X_{i})\rho_{0}(u/g(x)) w(x,u)
/\{g(x)^{2}\alpha_{x}(u)\}.$ By the central limit theorem for martingales,
one gets (see, e.g.,   \cite{Nielsen1995}, Proposition 1)
\begin{eqnarray*}
&\displaystyle(nb)^{1/2}\sum_{i=1}^{n}\int h_{ni}(x,u)\,\mathrm{d}M_{i}(u)
\Longrightarrow
N(0,\sigma
^{2}),&
\\
&\hspace*{-158pt}\displaystyle\sigma^{2} =
p\lim_{n\rightarrow\infty}
nb\sum_{i=1}^{n}\int h_{ni}^{2}(u)\,\mathrm{d} \langle M_{i}(u) \rangle&\\
&\hspace*{20pt}\displaystyle=
p\lim_{n\rightarrow\infty}
n^{-1}b\sum_{i=1}^{n}K_{b}^{2}(x-X_{i})\int\frac{\rho
_{0}^{2}(u/g(x))}{%
\alpha_{x}(u)^{2}g(x)^{4}} w(x,u) \alpha_{X_{i}}(u)C_{i}(u)Z_{i}(u)\,\mathrm{d}u&
\\
&\hspace*{-124pt}\displaystyle=
\|K\|_{2}^{2}\int\frac{\rho_{0}^{2}(u/g(x))}{\alpha
_{x}(u)g(x)^{4}}%
\varphi_{x}(u) w(x,u) \,\mathrm{d}u.&
\end{eqnarray*}

Furthermore,%
\begin{eqnarray*}
&&\int\rho_{0} \biggl\{ \frac{y}{g(x)} \biggr\} \frac{1}{g(x)^{2}}\frac{%
\mathcal{B}_{x,y}}{\alpha_{x}(y)} w(x,y) \,\mathrm{d}y \\
&&\quad =\frac{1}{n}\sum_{i=1}^{n}K_{b}(x-X_{i})\int\!\!\int\rho_{0} \biggl\{
\frac{y
}{g(x)} \biggr\}
\frac{w(x,y)}{\alpha_{x}(y)g(x)^{2}}k_{h}(y-y^{\prime})\\
&& \hphantom{\quad =\frac{1}{n}\sum_{i=1}^{n}K_{b}(x-X_{i})\int\!\!\int}
 {}\times
[ \alpha_{X_{i}}(y^{\prime})-\alpha_{x}(y) ] C_{i}(y^{\prime
})Z_{i}(y^{\prime})\,\mathrm{d}y^{\prime}\,\mathrm{d}y
\end{eqnarray*}
and it is easily seen that this can be written as a bias term of order $
\mathrm{O}(h^2) + \mathrm{O}(b^2)$, plus a remainder term of order $\mathrm{o}_P((nb)^{-1/2})$.

Finally, note that for any sequence $\delta_{n}\rightarrow0$, we have
%
%eB.8 ###
\begin{equation}
\sup_{|\theta-\theta_{0}|\leq\delta_{n}} \vert\hat
{l}_{\alpha
_{0}}^{\prime\prime}(\theta;x )-l^{\prime\prime}(\theta
_{0};x) \vert=\mathrm{o}_{P}(1), \label{l2}
\end{equation}
where%
\[
l^{\prime\prime}(\theta_{0};x)=2\int\rho_{0}^{2} \biggl\{ \frac{u}{g(x)}
\biggr\} \frac{w(x,u)}{g(x)^{4}}\frac{\varphi_{x}(u)}{\alpha_{x}(u)}\,\mathrm{d}u
\]
and $l^{\prime\prime}(\theta_{0};x)>0.$

From (\ref{te}), we then obtain
\[
\hat{g}_{C}^{o}(x)-g(x)=- \{ l^{\prime\prime}(\theta
_{0};x) \} ^{-1}\hat{l}_{\alpha_{0}}^{\prime}(\theta
_{0};x) \{ 1+\mathrm{o}_{P}(1) \}
\]
and the asymptotic distribution follows.
\end{pf*}
%
%EndExpansion

\begin{pf*}{Proof of Theorem \ref{alpha0}} We first consider the
infeasible estimator $\hat{\alpha}_{0,A}^{o}(y)$. Consider the following
decomposition:
\begin{eqnarray*}
\hat{\alpha}_{0,A}^{o}(y)
&=&\frac{\int E^A_{x_{,}yg(x)} w(x,yg(x))
\,\mathrm{d}x}{\int
 ({E^A_{x,yg(x)} w(x,yg(x))}/({g(x)\hat{\alpha}_{x,A} ( yg(x) ) }))
\,\mathrm{d}x}\\
&=&\frac{\hat{A}^{o}(y)}{\hat{B}^{o}(y)}\\
 &=& \frac{\hat
{A}^{o}(y)}{\int
 {%
\hat{B}^{o}_{1x}(y)}/{\hat{B}^{o}_{2x}(y)} \,\mathrm{d}x}\qquad \mbox{(say)} \\
&=&
\biggl[\frac{\hat{A}^{o}(y)}{B^{o}(y)}-\frac{A^{o}(y)}{B^{o}(y)^2}\int
\frac{\hat{B}_{1x}^{o}(y)}{B_{2x}^{o}(y)}\,\mathrm{d}x+\frac
{A^{o}(y)}{B^{o}(y)^2}\int
\frac{B_{1x}^{o}(y)\hat{B}_{2x}^{o}(y)}{B_{2x}^{o}(y)^2}\,\mathrm{d}x \biggr]
\bigl(1+\mathrm{o}_P(1)\bigr),
\end{eqnarray*}
where $A^o(y) = E[(CZ)\{yg(X)\} w\{X,yg(X)\}]$, $B^o(y) = E[(CZ)\{
yg(X)\}
w\{X,yg(X)\}] / \alpha_0(y)$,$\linebreak$ $B^o_{1x}(y) = E[(CZ)\{
yg(x)\}
w\{x,yg(x)\}|X=x]$ and $B^o_{2x}(y) = \alpha_0(y)$ are the limits of the
corresponding quantities with hats.

Straightforward calculations show that\vspace*{-2pt}
\begin{eqnarray*}
\hat{A}^{o}(y) &=& n^{-1}\sum_{i=1}^{n} ( C_{i}Z_{i} ) \{
yg(X_{i}) \} w\{X_i,yg(X_i)\} + \mathrm{o}_P((nh)^{-1/2}) + \mathrm{O}(h^{2}) + \mathrm{O}(b^2)
\end{eqnarray*}
and\vspace*{-2pt}
\begin{eqnarray*}
\int\frac{\hat{B}_{1x}^{o}(y)}{B_{2x}^{o}(y)}\,\mathrm{d}x
&=& \alpha_0(y)^{-1}
n^{-1}\sum_{i=1}^{n} ( C_{i}Z_{i} ) \{ yg(X_{i}) \}
w\{X_i,yg(X_i)\}\\
&&{} + \mathrm{o}_P((nh)^{-1/2})
 + \mathrm{O}(h^{2}) + \mathrm{O}(b^2).
\end{eqnarray*}
Next, we consider the term $\int B_{1x}^{o}(y)\hat{B}_{2x}^{o}(y)
B_{2x}^{o}(y)^{-2} \,\mathrm{d}x$. Decomposing $\hat\alpha_{x,A}(yg(x)) =
O^A_{x,yg(x)}/ E^A_{x,yg(x)}$ in $\hat{B}_{2x}^{o}(y)$ in a similar way as
above, we obtain, after some calculations, that\vspace*{-1pt}
\begin{eqnarray*}
&& \int\frac{B_{1x}^{o}(y)\hat{B}_{2x}^{o}(y)}{B_{2x}^{o}(y)^2} \,\mathrm{d}x \\[-1pt]
&&\quad  = B^o(y) + \frac{1}{\alpha_0(y)} n^{-1} \sum_{i=1}^n k_h\{
yg(X_i)-Y_i\}
C_i(Y_i) w(X_i,yg(X_i)) \\[-1pt]
&&\qquad{} - \frac{1}{\alpha_0(y)} n^{-1} \sum_{i=1}^n
(C_iZ_i)\{yg(X_i)\} w\{X_i,yg(X_i)\} + \mathrm{o}_P((nh)^{-1/2}) + \mathrm{O}(h^2) + \mathrm{O}(b^2).
\end{eqnarray*}
Putting the three terms together, we get that\vspace*{-1pt}%\vspace*{-9pt}
\begin{eqnarray*}
\hat\alpha_{0,A}^o(y)
 & = & \frac{1}{B^{o}(y)} n^{-1} \sum_{i=1}^n
k_h\{yg(X_i)-Y_i\} C_i(Y_i) w(X_i,yg(X_i)) \\
&&{} + \mathrm{o}_P((nh)^{-1/2}) + \mathrm{O}(h^2) + \mathrm{O}(b^2).
\end{eqnarray*}
We now consider the feasible estimator $\hat\alpha_{0,A}(y)$. Write\vspace*{-1pt}
\begin{eqnarray*}
 \hat{\alpha}_{0,A}^{o}(y)-E\hat{\alpha}_{0,A}^{o}(y)
 &=& \frac{1}{B^{o}(y)}\int k_{h} \{ yg(u)-v \} w(u,yg(u)) \,\mathrm{d} \bigl(
\hat{F}^o(u,v)-F^o(u,v) \bigr) \\
  &=& \frac{N^{o}(y)}{B^{o}(y)}\qquad \mbox{(say)},
\end{eqnarray*}
where $\hat F^o(u,v) = n^{-1} \sum_{i=1}^n I(X_i \le u, Y_i \le v,
C_i(Y_i)=1)$ and $F^o(u,v) = E[\hat F^o(u,v)]$. Therefore,
%
%eB.9 ###
\begin{eqnarray} \label{alphan}
&& [ \hat{\alpha}_{0,A}(y)-E\hat{\alpha}_{0,A}(y) ] - [ \hat{%
\alpha}_{0,A}^{o}(y)-E\hat{\alpha}_{0,A}^{o}(y) ] \nonumber
\\[-8pt]
\\[-8pt]
&&\quad  = \biggl[ \frac{1}{B(y)}-\frac{1}{B^{o}(y)} \biggr] N(y)+\frac{1}{B^{o}(y)}%
[ N(y)-N^{o}(y) ],\nonumber
\end{eqnarray}
where $B(y)$ and $N(y)$ are defined by replacing $g(\cdot)$ in the formulas
of $B^o(y)$ and $N^o(y)$ by $\hat g(\cdot)$, and where $E\hat{\alpha}
_{0,A}(y)$ is the expected value of $\hat{\alpha}_{0,A}(y)$ with
$\hat g$
considered as fixed. Write
\[
N^{o}(y)=h^{-1}\int k(z) w(u,v+hz)\, \mathrm{d} \bigl( [\hat{F}^o-F^o] \bigl(
u,yg(u)-hz \bigr) \bigr).
\]
Hence,
\begin{eqnarray*}
N(y)-N^{o}(y)
 &=&\mathrm{O}\Bigl ( h^{-1}\sup_{u, \vert t_{2}-t_{1} \vert
\leq C(nb)^{-1/2} (\log n)^{1/2}} \vert
\hat{F}^o(u,t_{1})-F^o(u,t_{1})\\
&&\hphantom{\mathrm{O}\Bigl ( h^{-1}\sup_{u, \vert t_{2}-t_{1} \vert
\leq C(nb)^{-1/2} (\log n)^{1/2}} \vert}{}-
\hat{F}^o(u,t_{2})+F^o(u,t_{2}) \vert\Bigr) \\
&=&\mathrm{O}_{P} ( h^{-1}n^{-1/2} (nb)^{-1/4} (\log n)^{1/4} ) =\mathrm{o}_{P} ( (nh)^{-1/2} ) ,
\end{eqnarray*}
provided $nh^{2}b (\log n)^{-1} \rightarrow\infty$. Next, note that
$N(y) =
\mathrm{O}_P((nh)^{-1/2})$, $B(y)-B^o(y) = \mathrm{O}_P((nb)^{-1/2}) = \mathrm{o}_P(1)$ and
hence (\ref{alphan}) is $\mathrm{o}_P((nh)^{-1/2})$. Since it can be easily
seen that
$E\hat{\alpha}_{0,A}(y)-E\hat{\alpha}_{0,A}^{o}(y) = \mathrm{O}(b^2) +
\mathrm{O}(h^2)$, it
follows that $\hat\alpha_{0,A}(y)$ and $\hat\alpha_{0,A}^o(y)$ are
asymptotically equivalent.

Finally, we consider the calculation of the asymptotic variance of
$\hat
\alpha_{0,A}(y)$:
\begin{eqnarray*}
  \operatorname{AsVar}(\hat\alpha_{0,A}(y))  &=& \frac{n^{-1}}{B^o(y)^2} \mbox{Var}\lbrack k_h\{yg(X)-Y\} C(Y)
w(X,yg(X))] \\
& =& \frac{n^{-1}}{B^o(y)^2} \int\!\!\int k_h^2\{yg(x)-t\} E[C(t)|Y=t,X=x]
\,\mathrm{d}F_x(t)\\
&&\hphantom{\frac{n^{-1}}{B^o(y)^2} \int\!\!\int}{}\times w^2(x,yg(x)) \,\mathrm{d}F(x)\bigl (1+\mathrm{o}(1)\bigr) \\
& =& \frac{(nh)^{-1}}{B^o(y)^2} \int k^2(u) \,\mathrm{d}u E \{
E[C(yg(X))|Y=yg(X),X]\\
&&\hphantom{\frac{(nh)^{-1}}{B^o(y)^2} \int}{}\times f_X(yg(X)) w^2(X,yg(X)) \} \bigl(1+\mathrm{o}(1)\bigr).
\end{eqnarray*}
\upqed
\end{pf*}

\begin{pf*}{Proof of Theorem \ref{g2step}} Consistency of $\widehat{%
\theta}=\hat{g}^{2-\mathit{step}}(x)$ follows similarly as in the proof of
Theorem \ref{goracle} from condition (A1), (\ref{Q1}) and the fact
that%
%
%eB.10 ###
\begin{equation}
\sup_{\theta\in I_n(x) } \vert\hat{l}_{\hat{\alpha
}}(\theta;x
) - \hat{l}_{\alpha_{0}}(\theta;x) \vert\stackrel
{P}{\rightarrow
}%
0. \label{Q1add}
\end{equation}
Equation (\ref{Q1add}) follows from assumption (C1) and the uniform
consistency of $\hat{\alpha}_{x,C}(y);$ see the proof of Theorem
\ref%
{goracle}. For the proof of Theorem \ref{g2step}, it remains to show
for $%
\theta_0=g(x)$ that for some $\gamma^*$,
%
%eB.12 ###
%eB.11 ###
\begin{eqnarray} \label{claim1}
\hat{l}_{\hat\alpha}^{\prime}(\theta_0;x)
&=& \hat
{l}_{\alpha
_{0}}^{\prime}(\theta_0;x)+h^2 \gamma^* + \mathrm{o}_P((nb)^{-1/2}),
\\ \label{claim2}
\hat{l}_{\hat\alpha}^{\prime\prime}(\theta;x ) &=&
\hat
{l}%
_{\alpha_{0}}^{\prime\prime}(\theta;x ) +\mathrm{o}_P(1),
\end{eqnarray}
uniformly for $\theta$ in a neighborhood of $\theta_0$. Claim (\ref{claim2})
follows immediately from assumption (C1). For the proof of (\ref{claim1}),
note first that
\begin{eqnarray*}
&&\hat{l}_{\hat\alpha}^{\prime}(\theta_0;x)- \hat
{l}_{\alpha
_{0}}^{\prime}(\theta_0;x) \\
&&\quad = 2\int\biggl[ \hat{\alpha}_{x,C}(y)-{\frac{1 }{\theta_0}} \alpha
_0 \biggl({\frac{y }{\theta_0}} \biggr) \biggr] (\widehat\rho_{0}- \rho
_{0})\biggl ( \frac{y}{\theta_0} \biggr) \frac{1}{\theta_0^{2}}\frac
{E_{x,y}^{C}%
}{\hat{\alpha}_{x,C}(y)} w(x,y) \,\mathrm{d}y \\
&&\qquad{}+ 2\int\biggl[ {\frac{1 }{\theta_0}}(\hat\alpha_0 -\alpha_{0}) \biggl({%
\frac{y }{\theta_0}} \biggr) \biggr] (\widehat\rho_{0}- \rho_{0}) \biggl(
\frac{y}{\theta_0} \biggr) \frac{1}{\theta_0^{2}}\frac
{E_{x,y}^{C}}{\hat{%
\alpha}_{x,C}(y)} w(x,y) \,\mathrm{d}y \\
&&\qquad{}+ 2\int\biggl[ {\frac{1 }{\theta_0}}(\hat\alpha_0 -\alpha_{0}) \biggl({%
\frac{y }{\theta_0}} \biggr) \biggr] \rho_{0} \biggl( \frac{y}{\theta_0} \biggr)
\frac{1}{\theta_0^{2}}\frac{E_{x,y}^{C}}{\hat{\alpha
}_{x,C}(y)} w(x,y)
\,\mathrm{d}y ,
\end{eqnarray*}
where $\widehat\rho_{0}(u)=\hat\alpha_{0}(u)+u\hat\alpha
_{0}^{\prime}(u)$. It follows from (C1) that the second term of the
right-hand side is of order $\mathrm{o}_P(n^{-2/5})$. From (C1) and (C2), we get
that up to
a deterministic term of order $\mathrm{O}(h^2)$, the third term is also of order
$%
\mathrm{o}_P(n^{-2/5})$. The first term is equal to $T_n + \mathrm{o}_P(n^{-2/5})$, where
\[
T_n = 2\int\biggl[ \hat{\alpha}_{x,C}(y)-{\frac{1 }{\theta_0}}
\alpha
_0 \biggl({\frac{y }{\theta_0}} \biggr) \biggr] (\widehat\rho_{0}- \rho
_{0}) \biggl( \frac{y}{\theta_0} \biggr) \frac{1}{\theta_0^{2}}\frac
{E_{x,y}^{C}%
}{{\alpha}_{x}(y)} w(x,y) \,\mathrm{d}y.
\]
For the proof of Theorem \ref{g2step}, it remains to show that
%
%eB.13 ###
\begin{equation}
T_n = \mathrm{o}_P(n^{-2/5}). \label{claim3}
\end{equation}
By application of (\ref{NL}), we can write $T_n= T_{n,1} + T_{n,2}$, where
\begin{eqnarray*}
T_{n,1} &=& 2\int\mathcal{V}_{x,y} (\widehat\rho_{0}- \rho_{0}) \biggl(
\frac{y}{\theta_0} \biggr) \frac{1}{\theta_0^{2}}\frac{1}{{\alpha}_{x}(y)}
w(x,y) \,\mathrm{d}y, \\
T_{n,2} &=& 2\int\mathcal{B} _{x,y} (\widehat\rho_{0}- \rho_{0}) \biggl(
\frac{y}{\theta_0} \biggr) \frac{1}{\theta_0^{2}}\frac{1}{{\alpha}_{x}(y)}
w(x,y) \,\mathrm{d}y.
\end{eqnarray*}

It can be easily checked that $T_{n,2}=\mathrm{o}_{P}(n^{-2/5})$ (cf. the
proof of Theorem \ref{goracle}). The term $T_{n,1}$ can be decomposed
into $%
T_{n,11}+T_{n,12}$, where
\[
T_{n,11} =\sum_{i=1}^{n}\int h_{ni}(x,u)\,\mathrm{d}M_{i}(u),\qquad
T_{n,12} =\sum_{i=1}^{n}\int g_{ni}(x,u)\,\mathrm{d}M_{i}(u),
\]
with
\begin{eqnarray*}
h_{ni}(x,u) &=&\frac{2}{n}K_{b}(x-X_{i})\int\biggl[ \int(\hat{\alpha}
_{0}-\alpha_{0})\biggl \{ \frac{y}{\theta_{0}} \biggr\} \frac{1}{\theta
_{0}^{2}}\frac{1}{\alpha_{x}(y)}k_{h}(y-u)w(x,y)\,\mathrm{d}y \biggr] , \\
g_{ni}(x,u) &=&\frac{2}{n}K_{b}(x-X_{i})\int\biggl[ \int\biggl\{ \frac{y}{%
\theta_{0}} \biggr\} (\hat{\alpha}_{0}^{\prime}-\alpha_{0}^{\prime
}) \biggl\{ \frac{y}{\theta_{0}} \biggr\} \frac{1}{\theta_{0}^{2}}\frac{1}{%
\alpha_{x}(y)}k_{h}(y-u)w(x,y)\,\mathrm{d}y \biggr] .
\end{eqnarray*}
We now show that $T_{n,12}=\mathrm{o}_{P}(n^{-2/5})$. The claim $%
T_{n,11}=\mathrm{o}_{P}(n^{-2/5})$ can be shown by similar methods. For the
proof, we
apply   \cite{Geer2000}, Lemma 5.14. This lemma gives a bound on the
increments of the empirical process applied to function classes that depend
on the sample size. We apply the lemma with a fixed value of $x$, conditional
on the event that the number of values of $X_{i}$ in the support of $K_{b}$
is equal to $m$, where $m$ is of the same order as $nb$. We consider the
class of functions $g\dvtx J(x)\rightarrow\mathbf{R}$ such that, with a
sufficiently large constant $C$, for all $z\in J(x)$, $|g(z)-\alpha
_{0}^{\prime}(z)|\leq C\delta_{1,n}$ and $|g^{\prime}(z)|\leq
C\delta
_{2,n}$. We apply the lemma with $\alpha=\beta=1$ and $M=C\delta_{2,n}$.
We get that
\[
\sup\Biggl\vert\frac{1}{n}\sum_{i=1}^{n}K_{b}(x-X_{i})\int\biggl[ \int
\biggl( \frac{y}{\theta_{0}} \biggr) g \biggl( \frac{y}{\theta_{0}} \biggr)
\frac{1}{\theta_{0}^{2}}\frac{1}{\alpha_{x}(y)}k_{h}(y-u)w(x,y)\,\mathrm{d}y \biggr]
\,\mathrm{d}M_{i}(u) \Biggr\vert
\]
is of order $\mathrm{O}_{P}(\delta_{1,n}^{1/2}\delta
_{2,n}^{1/2}(nb)^{-1/2}+\delta
_{2,n}(nb)^{-1})=\mathrm{o}_{P}((nb)^{-1/2})$. This shows that $T_{n,12}=\mathrm{o}_{P}(n^{-2/5})$
and thus concludes the proof of Theorem \ref{g2step}.
\end{pf*}
\end{appendix}

\section*{Acknowledgements}
Research of O. Linton was supported by the ESRC.
Research of E. Mammen was supported by the Deutsche
Forschungsgemeinschaft, Project MA1026/11-1.
Research of I. Van Keilegom was supported by IAP research network grant
no. P6/03 of the Belgian government
(Belgian Science Policy) and by the European Research Council under the
European Community's Seventh Framework Programme
(FP7/2007-2013)/ERC Grant agreement no. 203650. The authors would
like to thank Sorawoot Srisuma for research assistance.

\printhistory

\end{document}